\documentclass{aims} % Use the aims.cls file to compile your paper
\usepackage[pagewise]{lineno}%\linenumbers
\usepackage{amsmath}
\usepackage{paralist}
\usepackage[misc]{ifsym}
\usepackage{epsfig} %For pictures: screened artwork should be set up with an 85 or 100 line screen
\usepackage{epstopdf} %This is to transfer .eps figure to .pdf figure; please compile your article using PDFLaTex or PDFTeXify.
\usepackage[colorlinks=true]{hyperref}
\hypersetup{urlcolor=blue, citecolor=red}
\usepackage{hyperref}
\usepackage{subfigure}
\usepackage[dvipsnames]{xcolor}

\allowdisplaybreaks

\def\currentvolume{X}
 \def\currentissue{X}
  \def\currentyear{200X}
   \def\currentmonth{XX}
    \def\ppages{X-XX}
     \def\DOI{10.3934/xx.xxxxxxx}

\usepackage{stmaryrd}
\usepackage{algorithm}
\usepackage{algorithmic}
\usepackage{multirow}
\usepackage{mathrsfs}
\usepackage{booktabs}
\newtheorem{theorem}{Theorem}[section]

\newtheorem{lemma}[theorem]{Lemma}
\newtheorem{proposition}[theorem]{Proposition}

\theoremstyle{definition}

\title[Super-resolution via low rank triple decomposition]
{Hyperspectral super-resolution via low rank tensor triple decomposition} % Only the first word and proper nouns should be capitalized.

\author[Xiaofei Cui and Jingya Chang]{}


\thanks{This work is supported by National Natural Science Foundation of China (grant No. 11901118 and 62073087).}

\thanks{$^*$Corresponding author: Jingya Chang}

\begin{document}
\maketitle

% Enter the first author's name and email address; email addresses are required for each author.
% Use footnote notations to indicate address and affiliations with commas between numbers if more than one address applies; see below for examples.
\centerline{\scshape
Xiaofei Cui$^{{\href{mailto:cuixiaofei3583@163.com}{\textrm{\Letter}}}1}$
and Jingya Chang$^{{\href{mailto:jychang@gdut.edu.cn.com}{\textrm{\Letter}}}*1}$}

\medskip

{\footnotesize
% Enter the full affiliation and country name:
% Do not insert commas or periods at the end of lines.
 \centerline{$^1$School of Mathematics and Statistics, Guangdong University of Technology, China}
} % Do not forget to end {\footnotesize with the sign }

%\medskip
%
%{\footnotesize
% % Enter the full affiliation and country name:
% \centerline{$^2$School of Mathematics and Statistics, Guangdong University of Technology, China}
%}

\bigskip

% The name of the handling editor will be entered by AIMS production staff.
% "Communicated by Handling Editor" is not needed for special issue.
% \centerline{(Communicated by Handling Editor)} 

%%%%%%%%%%%%%%%%%%%%%%%%%%%%%%%%%%%%%%%%%%%%%%%%%%%%%%%
%             5. ABSTRACT
%%%%%%%%%%%%%%%%%%%%%%%%%%%%%%%%%%%%%%%%%%%%%%%%%%%%%%%

\begin{abstract}
Hyperspectral image (HSI) and multispectral image (MSI) fusion aims at producing a super-resolution image (SRI). In this paper, we  establish  a nonconvex optimization model for image fusion problem through low rank tensor triple decomposition. Using the L-BFGS approach, we develop a first-order optimization algorithm for obtaining the desired super-resolution image (TTDSR). Furthermore, two detailed methods are provided for calculating the gradient of the objective function. With the aid of the Kurdyka-$\L$ojasiewicz property, the iterative sequence is proved to converge to a stationary point. Finally, experimental results on different datasets show the effectiveness of our proposed approach.
\end{abstract}

%%%%%%%%%%%%%%%%%%%%%%%%%%%%%%%%%%%%%%%%%%%%%%%%%%%%%%
%                   6. BODY
%%%%%%%%%%%%%%%%%%%%%%%%%%%%%%%%%%%%%%%%%%%%%%%%%%%%%%

% Only the first word and proper nouns of section titles should be capitalized.
% The title of section 1:
\section{Introduction}

Equipment for hyperspectral imaging uses an unusual sensor to produce images with great spectral resolution. The key factor in an image is the electromagnetic spectrum obtained by sampling at hundreds of continuous wavelength intervals. Because hyperspectral images have hundreds of bands, spectral information is abundant \cite{32}. Hyperspectral images are now the core of a large and increasing number of remote sensing applications, such as environmental monitoring \cite{4,25}, target classification \cite{41,8} and anomaly detection \cite{7,37}.
	
With optical sensors, there is a fundamental compromise among spectral resolution, spatial resolution and signal-to-noise ratio. Hyperspectral images therefore have a high spectral resolution but a low spatial resolution. Contrarily, multispectral sensors, which only have a few spectral bands, can produce images with higher spatial resolution but poorer spectral resolution \cite{5}.

The approaches for addressing hyperspectral and multispectral images fusion issues can be categorized into four classes \cite{42}. The  pan-sharpening approach falls within the first category, which includes component substitution \cite{27,1} and multi-resolution analysis \cite{6}. Selva \cite{30} and Liu \cite{22} proposed multi-resolution analysis method to extract spatial details that were injected into the interpolated hyperspectral image. Chen et al. \cite{9} divided all bands of the images into several regions and then fused the images of each region by the pan-sharpening approach. Unfortunately, these methods may lead to spectral distortion. The second category is subspace-based formulas. The main idea is to explore the low-dimensional representation of hyperspectral images, which are regarded as a linear combination of a set of basis vectors or spectral features \cite{18,38,34}. Hardie et al. \cite{15} employed Bayesian formulas to solve the image fusion problem with the maximum posteriori probability. The method based on matrix factorization falls within the third category. Under these circumstances, HSI and MSI are regarded as the degradation versions of the reference SRI. Yokoya et al. \cite{40} proposed a coupled nonnegative matrix factorization method, which estimated the basic vectors and its coefficients blindly and then calculated the endmember matrix and degree matrix by alternating optimization algorithm. In addition, there are methods combining  sparse regularization and spatial regularization \cite{19,31}, such as a low rank matrix decomposition model for decomposing a hyperspectral image into a super-resolution image and a difference image in \cite{12,42}.	

Tensor decomposition methods fall within the fourth category. Tensor, a higher-order extension of matrix, can capture the correlation between the spatial and spectral domains in hyperspectral images simultaneously. Representing hyperspectral and multispectral images as third-order tensors has been widely used in hyperspectral image denoising and super-resolution problems. The hyperspectral image is a three-dimensional data cube with two spatial dimensions (height and width) and one spectral dimension \cite{17,33,23}. The spatial details and spectral structures of the SRI therefore are both more accurate.
	
By utilizing low rank tensor decomposition technology to explore the low rank properties of hyperspectral images, high spatial resolution and accurate restoration guarantees are obtained. Kanatsoulis et al. \cite{20} proposed treating the super-resolution problem as a coupled tensor approximation, where the super-resolution image satisfied low rank Canonical Polynomial (CP) decomposition. Later, Li et al. \cite{28,24} considered a novel approach based on tucker decomposition. Xu et al. \cite{39} first considered that the target SRI exhibited both the sparsity and the piecewise smoothness on three modes and then used tucker decomposition to construct model. Dian et al. \cite{13} applied tensor kernel norm to explore low rank property of SRI. He et al. \cite{16} proposed tensor ring decomposition model for addressing the image fusion problem.

\begin{figure}[h]
		\centering
		\begin{minipage}[t]{0.45\textwidth}
			\centering
			\includegraphics[width=3.5cm]{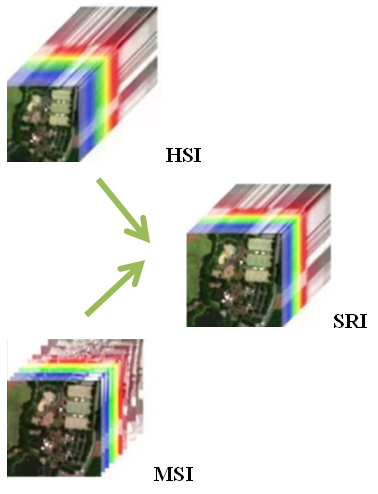}
		\end{minipage}
		\begin{minipage}[t]{0.45\textwidth}
			\centering
			\includegraphics[width=6.2cm]{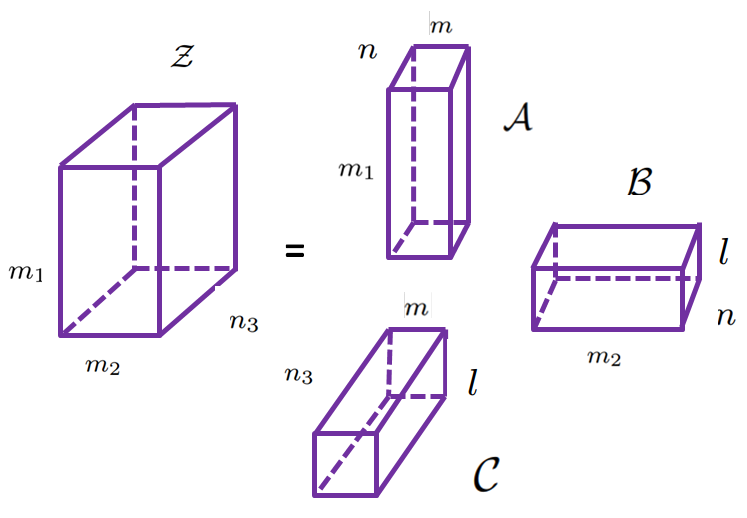}
		\end{minipage}
		\caption{The image fusion is displayed in the left image, and the tensor triple decomposition structure is shown in the right.}\label{fig.1}
	\end{figure}

In this paper, we propose a novel model based on low rank tensor triple decomposition to address the hyperspectral and multispectral images fusion problem. Triple decomposition breaks down a third-order tensor into three smaller tensors with predetermined dimensions. It is effectively applied in tensor completion and image restoration \cite{29,11}. The structures of image fusion and triple decomposition  are presented in Figure \ref{fig.1}. The optimization model is solved by the L-BFGS algorithm. The convergent results are given with the help of Kurdyka-$\L$ojasiewicz property of the objective function. Numerical experiments demonstrate that the  proposed TTDSR method performs well and the convergent conclusion is also validated numerically.
	
The rest of this paper is organized as follows. In section \ref{section:2}, we first provide definitions used in the paper and present the basic algebraic operations of tensors. In section \ref{section:3}, we introduce  our motivation and model for the hyperspectral and multispectral images fusion. In section \ref{section:4}, the L-BFGS algorithm and the gradient calculation are also presented in detail. Convergence analysis of the algorithm is available in section \ref{section:5}. In section \ref{section:6}, experimental results on two datasets are given. The conclusion is given in section \ref{section:7}.

% The title of section 2:
\section{Preliminary}
\label{section:2}
% The title of the first subsection in section 2:
Let $\mathbb{R}$ be the real field. We use lowercase letters and boldface lowercase letters  to represent scalars and vectors respectively, while capital letters and calligraphic letters stand for matrices and tensors respectively.

\begin{figure}[htp]
\begin{center}
  % replace aims_logo.pdf by your figure file name
  \includegraphics[width=1\textwidth]{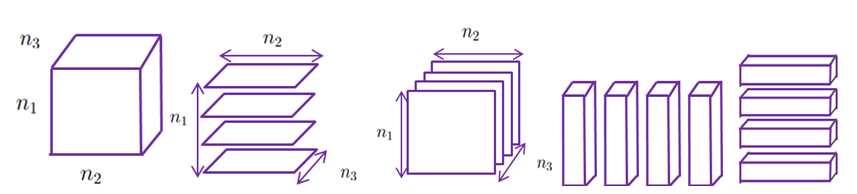}
  \caption{Slices and fibers of tensor $\mathcal{A}$. From left to right are horizontal slices $\mathcal{A}_{i,:,:}$, frontal slices $\mathcal{A}_{:,:,k}$, column fibers $\mathcal{A}_{i,:,k}$, and row fibers $\mathcal{A}_{:,j,k}$.}\label{fig.2}
  \end{center}
\end{figure}

	An $m$th order $n$ dimensional tensor $\mathcal{A}$ is a multidimensional array with its entry being
\[
\mathcal{A}_{i_1,i_2,\ldots,i_m}=a_{i_1,i_2,\ldots,i_m},\ \text{for} \ i_j=1,\ldots,n, j=1,\ldots,m.
\]
Similarly, the $i$th entry in a vector $\textbf{a}$ is symbolized by $(\mathbf{a})_{i}=a_{i}$, and the $(i,j)$th entry of a matrix  $A$  is  $(A)_{i j}=a_{i j}.$ Unless otherwise specified,  the order of tensor  is set to three  hereafter in this paper. Without loss of generality, we denote $\mathcal{A}\in \mathbb{R}^{n_1 \times n_2 \times n_3}.$ A fiber of  $\mathcal{A} $  is a vector produced by  fixing two indices of $\mathcal{A},$ such as $\mathcal{A}(i,j,:)$ for any $i\in[1,\ldots,n_1],j\in[1,\ldots,n_2]$. A slice of   $\mathcal{A}$  is a matrix by varying two of its indices while  fixing another one, such as $\mathcal{A}(:,:,k)$ for any $k\in[1,\ldots,n_3].$
The geometric structures of the tensor expansion are shown in Figure \ref{fig.2}.

The mode-$k$ product of the tensor $\mathcal{A}$ and a matrix $X$ is an extension of matrix product.
Suppose  matrices $ F \in \mathbb{R}^{l \times n_1}$, $G \in \mathbb{R}^{m\times n_2}$, and  $H \in \mathbb{R}^{n \times n_3} $. The mode-$1$ product of $\mathcal{A}$ and $F$ is denoted as $\mathcal{A}\times_{1} F$ with its elements being
$$\left(\mathcal{A} \times_{1} {F}\right)_{t j k}=\sum_{i=1}^{n_1} a_{i j k} f_{t i}, \quad \text{for} \ t = 1,\ldots l, j=1,\ldots,n_2, k=1,\ldots n_3. $$
Also we have
	\begin{equation}
	\left(\mathcal{A} \times_{2} {G}\right)_{i p k}=\sum_{j=1}^{n_2} a_{i j k} g_{p j},
    \quad\left(\mathcal{A} \times_{3} {H}\right)_{i j q}=\sum_{k=1}^{n_3} a_{i j k} h_{q k},\quad
	\end{equation}
\text{for} $ p = 1,\ldots m,$ $q=1,\ldots,n.$ It is easy to verify that
	\begin{equation}
\begin{split}
  \mathcal{A} \times_{p} {F} \times_{q} {G}= & \mathcal{A}\times_{q} {G} \times_{p} {F}, \quad \forall \ p,q=1,2,3 \ \text{and} \ p \neq q, \\
   \mathcal{A}\times_{k} {F} \times_{k} {G}= & \mathcal{A}\times_{k}({GF}), \quad \forall \ k = 1,2,3.
\end{split}
	\end{equation}
The mode-$k$  matricization of $ \mathcal{A} $ denoted by  ${A}_{(k)}$  arranges its $n_k$ mode-$k$  fibers  as the columns of $A_{(k)}$ in order.
The Frobenius norm of tensor $\mathcal{A}$ is given by
	\begin{equation}
	\|\mathcal{A}\|_{F}=\left(\sum_{i=1}^{n_1} \sum_{j=1}^{n_2} \sum_{k=1}^{n_3} a_{i j k}^{2}\right)^{1 / 2}.
	\end{equation}
	
	The unfolding related no matter with tensors, matrices or vectors in the following part is carried out under the precondition that the left index varies more rapidly than the right one. Subscript represents dimension. We note that $ {vec}\left(B_{m\times n }\right) $ creates a column vector  $\mathbf{b}_{m n} $ by stacking the columns of  $B_{m \times n} $ below one another, i.e.,  $\left(\mathbf{b}_{m n}\right)_{i+(j-1) m}=b_{i j}$.
	In addition, the kronecker product of matrices $ B_{m \times n}$  and  $C_{p \times q}$ is
	$$
	B_{m \times n} \otimes C_{p \times q}=\left(\begin{array}{cccc}
	b_{11} C_{p \times q} & b_{12} C_{p \times q} & \cdots & b_{1n} C_{p \times q}\\
	b_{21} C_{p \times q} & b_{22} C_{p \times q} & \cdots & b_{2n} C_{p \times q}\\
	\cdots & \cdots & \cdots & \cdots\\
	b_{m1} C_{p \times q} & b_{m2} C_{p\times q} & \cdots & b_{mn} C_{p \times q}
	\end{array}\right), \quad
	$$
\text{for} $ s=1,\cdots,p,$ $t=1,\cdots,q.$ It generates a large matrix $ D_{pm \times qn}=B_{m \times n} \otimes C_{p \times q} $, whose entries are $ {d}_{gh}=b_{i j} c_{s t} $ with $ g=s+(i-1) p $ and $ h=t+(j-1) q $. By this rule, unfolding matrices of tensor $\mathcal{A}$ could be represented by
	\begin{equation}
	\begin{split}
	\left(A_{n_1 \times n_2 n_3}\right)_{i, j+(k-1) n_2}=a_{i j k},\quad &\left(A_{n_1 \times n_3 n_2}\right)_{i, k+(j-1) n_3}=a_{i j k}, \\
	\left(A_{n_2 \times n_3 n_1 }\right)_{j, k+(i-1) n_3}=a_{i j k}, \quad &\left(A_{n_2 \times n_1  n_3}\right)_{j, i+(k-1) n_1 }=a_{i j k}, \\
	\left(A_{n_3 \times n_1  n_2}\right)_{k, i+(j-1) n_1 }=a_{i j k},\quad &\left(A_{n_3 \times n_2 n_1 }\right)_{k, j+(i-1) n_2}=a_{i j k}.
	\end{split}
	\end{equation}
	Moreover, $A_{n_1 n_2 \times n_3}=\left(A_{n_3 \times n_1 n_2}\right)^{\top}.$

\section{The tensor triple decomposition (TTD) model for HSI-MSI fusion}
	\label{section:3}
Tensors can describe high dimensional relationship. The hyperspectral and multispectral images are naturally third-order tensors. Methods based on low rank matrix decomposition reconstruct the 3-d tensor data into a 2-d matrix data at beginning, which destroy the original 3-d structure and may cause negative influence on the fusion consequences. Therefore,  hyperspectral and multispectral image fusion problems have been investigated with the aid of  various tensor decompositions in recent years. The tensor triple decomposition model has certain advantages over the classical CP and tucker decompositions. The triple rank is less than or equal to CP rank, which has been proved theoretically and practically. As described in theorem 2.1 of \cite{29}, low rank triple decomposition and triple ranks  are well-defined. The proposed triple decomposition has the low rank properties and performs well in image restoration.  When facing large-scale problems, the lower rank indicates that we transform the problem into a lower dimensional subspace, which is relatively easy to solve. Thus, we employ the tensor triple decomposition framework to establish the hyperspectral super-resolution model.

Define $\mathcal{Y}_{h} \in \mathbb{R}^{n_{1} \times n_{2}\times n_{3}}$, $\mathcal{Y}_{m}\in \mathbb{R}^{m_{1} \times m_{2}\times m_{3}}$ and $\mathcal{Z} \in \mathbb{R}^{m_{1} \times m_{2}\times n_{3}}$ as  hyperspectral image (HSI), multispectral image (MSI) and super-resolution image (SRI) respectively. Here the first two indices $n_1$ and $n_2$ or $m_1$ and $m_2$ denote the spatial dimensions. The third index $n_3$ or $m_3$ indicate the number of spectral bands. Usually, the MSI contains more spatial information than the HSI, that is $m_1\gg n_1, m_2\gg n_2.$ The HSI has more spectral bands than MSI, which means $n_3\gg m_3.$ Our target is to find the SRI that possesses the  high spectral resolution of HSI and the spatial resolution of the MSI, i.e. a super-resolution image.

\subsection{The links among SRI, HSI and MSI}
Let the mode-$3$ matricization of $\mathcal{Y}_{h}$, $\mathcal{Y}_{m}$ and $\mathcal{Z}$ be $(Y_h)_{(3)}\in \mathbb{R}^{n_{1}n_2\times n_3 }, (Y_m)_{(3)}\in\mathbb{R}^{m_{1}m_2\times m_3} $ and $Z_{(3)}\in\mathbb{R}^{m_{1}m_2\times n_3}.$ The key point  which helps us to construct the relationship among SRI, HSI and MSI is that  there exist two linear operators $P_h\in\mathbb{R}^{n_{1}n_2\times m_1m_2 }$ and $P_3\in\mathbb{R}^{m_3\times n_3 }$ such that $(Y_h)_{(3)}=P_hZ_{(3)}$  and $(Y_m)_{(3)}=Z_{(3)}P_3^{T}$ \cite{20}. Thus we have
	\begin{equation}\label{3}
	\mathcal{Y}_{h}(:,:,k)=P_1\mathcal{Z}(:,:,k)P_2^T, \quad \text{for} \ k=1,\ldots, n_3,
	\end{equation}
	where $P_1\in\mathbb{R}^{n_1\times m_1}$, $P_2\in\mathbb{R}^{n_2\times m_2} $ and $P_2\otimes P_1=P_h.$ We also have
	\begin{equation}\label{4}
	\mathcal{Y}_{m}(i, j,:)={P}_{3} \mathcal{Z}(i, j,:), \quad  \text{for} \ i=1,\ldots, m_1, \ j=1,\ldots, m_2,
	\end{equation}
	where  $\mathcal{Z}(i, j,:) \in  \mathbb{R}^{n_{3}}$  represents a fiber of $\mathcal{Z}$ and  $\mathcal{Y}_{m}(i, j,:) \in \mathbb{R}^{m_{3}}$  is a fiber of  $\mathcal{Y}_{m}$ respectively. Therefore, $\mathcal{Y}_{h}$, $\mathcal{Y}_{m}$ can be rewritten as
	\begin{equation}\label{3.1}
	\mathcal{Y}_{h}=\mathcal{Z}\times_{1}{P}_{1}\times_{2}{P}_{2}, \quad\mathcal{Y}_{m}=\mathcal{Z}\times_{3}{P}_{3}.
	\end{equation}
	Equations $(\ref{3})$ and $(\ref{4})$ also are known as the spatial and spectral degradation  processes of $\mathcal{Z}.$  The matrices
	$P_{1}$, $P_{2}$ respectively describe the downsampling and blurring of the spatial degradation process. Downsampling is considered as linear compression, while blurring describes a linear mixing of neighbouring pixels under a specific kernel in both the row and column dimensions. The matrix ${P}_{3}$ is usually modeled as a band-selection matrix that selects the common spectral bands of the SRI and MSI.
	
	Based on the correlation of spatial and spectral information in hyperspectral and multispectral images, various low rank tensor decomposition models are established to study the problem. For example, suppose
	$\mathcal{Z} $ is decomposed by Canonical Polyadic decomposition into the sum of several rank-$1$ tensors, which is represented as
	$ \mathcal{Z}=\llbracket{A}, {B}, {C}\rrbracket.$ The HSI-MSI fusion model \cite{20} is denoted as
	$\mathcal{Y}_{h}=\llbracket{P}_{1} {A}, {P}_{2}{B}, {C}\rrbracket$ and $\mathcal{Y}_{m}=\llbracket{A}, {B}, {P}_{3} {C}\rrbracket$. Besides, tucker decomposition is considered \cite{28,44} and it is expressed as
	$\mathcal{Y}_{h}=\mathcal{C} \times_{1} \left({P}_{1} {A}\right)  \times_{2} \left({P}_{2} {B}\right) \times_{3} {C}$ and $\mathcal{Y}_{m}=\mathcal{C} \times_{1} {A} \times_{2} {B}  \times_{3} \left({P}_{3} {C}\right)$.

\subsection{The new model for HSI-MSI fusion}
\label{subsection:3.2}
	In \cite{29}, Qi et al. proposed a new tensor decomposition, named tensor triple decomposition, which can effectively express the original tensor information by three low rank tensors. Therefore, we establish the hyperspectral and multispectral image fusion model by low rank tensor triple decomposition.
	
	First, we introduce the tensor triple decomposition in detail. Tensor triple decomposition of a third-order $\mathcal{Z} \in \mathbb{Re}^{ m_{1}\times m_{2} \times n_{3}}$ representing the target SRI takes the following form
	\begin{equation}\label{2.7}
	\mathcal{Z}=\mathcal{A}\mathcal{B}\mathcal{C},
	\end{equation}
	where $\mathcal{A} \in \mathbb{Re}^{m_{1} \times m \times n}$, $\mathcal{B} \in \mathbb{Re}^{l \times m_{2} \times n } $ and  $\mathcal{C} \in \mathbb{Re}^{l \times m \times n_{3}}.$
	It can also be denoted as $ \mathcal{Z}=\llbracket \mathcal{A}, \mathcal{B}, \mathcal{C} \rrbracket$
	with its entries
	\begin{equation}\label{2.8}
	(\llbracket \mathcal{A}, \mathcal{B}, \mathcal{C} \rrbracket)_{i j k}=\sum_{t=1}^{l} \sum_{p=1}^{m} \sum_{q=1}^{n} a_{i p q} b_{t j q} c_{t p k},
	\end{equation}
	for $ i= 1, \ldots, m_{1}$, $j = 1, \ldots, m_{2}$, and $k = 1, \ldots, n_{3}.$
	The smallest value of  $r$ such that  (\ref{2.8}) holds is called the triple rank of  $ \mathcal{Z} $ and denoted  as  $\operatorname{rank}(\mathcal{Z})=r $. If $l=m=n=r$, equation (\ref{2.7}) is called low rank triple decomposition of $ \mathcal{Z}$, where $ r \leq \operatorname{mid}\left\{m_{1}, m_{2}, n_{3}\right\}$ and mid represents the larger number in the middle \cite{29,11}.	
	In particular, according to \cite[Theorem 2.2]{11}, triple decomposition satisfies the following equations
	\begin{equation}\label{3.10}
	\begin{split}
	Z_{m_1 \times m_2 n_3 }&=A_{m_1 \times n m}\left(E_{m} \otimes B_{n \times m_2 l}\right)\left(C_{l m \times n_3} \otimes E_{m_2}\right),\\
	Z_{m_2 \times n_3 m_1}&=B_{m_2 \times l n}\left(E_{n} \otimes C_{l \times n_3 m}\right)\left(A_{m n \times m_1} \otimes E_{n_3}\right),\\
	Z_{n_3  \times m_2 m_1}&=C_{n_3 \times l m}\left(E_{m} \otimes B_{l \times m_2  n}\right)\left(A_{n m \times m_1} \otimes E_{m_2 }\right),
	\end{split}
	\end{equation}
	where $E$  is an identity matrix with a proper size.
	
	Next, we propose a new model with the help of low rank tensor triple decomposition.
	Assume $\mathcal{A}$, $\mathcal{B} $ and $ \mathcal{C} $ are the low rank triple decomposition tensors of  tensor $ \mathcal{Z}.$ For the connection \eqref{3.1}, we have
	\begin{equation}
	\mathcal{Y}_{h}=\llbracket\mathcal{A}\times_{1}{P_{1}},\mathcal{B}\times_{2}{P_{2}},\mathcal{C}\rrbracket \quad \text{and} \quad %+ {N}_{h} \quad
\mathcal{Y}_{m}=\llbracket\mathcal{A},\mathcal{B},\mathcal{C}\times_{3}{P_{3}}\rrbracket. \label{3.2}  % + {N}_{m}.
	\end{equation}
    Since the best low rank approximation problem of tensors may be ill-posed, we add the Tikhonov regularization term and get
	the following optimization model
	\begin{equation}\label{3.3}
	\begin{split}
	\min_{\mathcal{A}, \mathcal{B}, \mathcal{C}} f(\mathcal{A}, \mathcal{B}, \mathcal{C}):&= \left\|\mathcal{Y}_{h}-\llbracket\mathcal{A}\times_{1}{P_{1}},\mathcal{B}\times_{2}{P_{2}},\mathcal{C}\rrbracket\right\|_{F}^{2}
	+\left\|\mathcal{Y}_{m}-\llbracket\mathcal{A},\mathcal{B},\mathcal{C}\times_{3}{P_{3}}\rrbracket\right\|_{F}^{2}\\
	&\quad+\mu\left(\|\mathcal{A}\|_{F}^{2}+\|\mathcal{B}\|_{F}^{2}+\|\mathcal{C}\|_{F}^{2}\right),
	\end{split}
	\end{equation}\\
	where $\mu$ is regularization parameter. Thus, we employ the optimization model  \eqref{3.3} to obtain the triple decomposition tensors $\mathcal{A},\mathcal{B}$ and $\mathcal{C}$, and to produce the  super-resolution image SRI by $ \mathcal{Z}=\llbracket \mathcal{A}, \mathcal{B}, \mathcal{C} \rrbracket$.

\section{The L-BFGS algorithm for solving the model}
	\label{section:4}
	In this section, we focus on numerical approaches for computing a first-order stationary point of the optimization problem \eqref{3.3}.
	Since the limited memory BFGS (L-BFGS) algorithm \cite{21} is powerful for large scale nonlinear unconstrained optimization, we apply it  to produce a search direction and consider using inexact line search techniques to update the iteration step size. In the computaion process, we either matricize or vectorize the tensor variable to get the gradient of the objective function. For convenience, we demonstrate the algorithm and its convergence analysis for
	 \begin{equation}
	\mathbf{x}:=\left(\begin{array}{l}
	\mathbf{a}_{m_{1} r r} \\
	\mathbf{b}_{m_{2} r r} \\
	\mathbf{c}_{n_{3} r r}
	\end{array}\right)=\left(\begin{array}{c}
	A_{m_{1} r r \times 1} \\
	B_{m_{2} r r \times 1} \\
	C_{n_{3} r r\times 1}
	\end{array}\right).
	\end{equation}
\subsection{Limited memory BFGS algorithm}
	BFGS is a quasi-Newton method which updates the approximation of the inverse of a Hessian iteratively. In the current iteration $k$, it constructs an approximation matrix $H_k$ to estimate the inverse of Hessian of $f(\mathbf{x})$.  The gradient of $f(\mathbf{x})$ is defined as $\mathbf{g(x)}.$  At the beginning, we introduce the basic BFGS update. Let
	\begin{equation}
	\quad \mathbf{y}_{k}=\nabla f\left(\mathbf{x}_{k+1}\right)-\nabla f\left(\mathbf{x}_{k}\right), \quad \mathbf{s}_{k}=\mathbf{x}_{k+1}-\mathbf{x}_{k}, \quad V_{k}=I-\rho_{k} \mathbf{y}_{k} \mathbf{s}_{k}^{\top},\label{4.1}
	\end{equation}
	and
	\begin{equation}\label{4.2}
	\begin{aligned}
	\begin{split}
	\rho_{k}=\left\{\begin{array}{ll}
	\frac{1}{\mathbf{y}_{k}^{\top} \mathbf{s}_{k}}, & \text { if } \mathbf{y}_{k}^{\top} \mathbf{s}_{k} \geq \epsilon\\
	0, & \text { otherwise }
	\end{array}\right.
	\end{split}
	\end{aligned},
	\end{equation}
	where  $I$  is an identity matrix, $\top$ represents transposition and  $\epsilon \in(0,1) $ is a small positive constant.
	The  matrix  $H_k$ is updated by
	\begin{equation}
	H_{k+1}=V_{k}^{\top} H_{k} V_{k}+\rho_{k} \mathbf{s}_{k} \mathbf{s}_{k}^{\top}.\label{4.3}
	\end{equation}
	To deal with large-scale optimization problems, Liu and Nocedal \cite{21} proposed the L-BFGS algorithm that implemented BFGS in an economical manner. Given a constant $l.$ When $k\geq l,$ only information of $l$ matrices $H_{k-1}, H_{k-2},\ldots,H_{k-l}$ are used to calculate the matrix $H_k$ by the following recursive form
		\begin{equation}
		{H_k}=V_{k-m}^{\top} H_{k-m} V_{k-m}+\rho_{k-m} \mathbf{s}_{k-m} \mathbf{s}_{k-m}^{\top}, \quad \text { for } m=1, 2, \ldots, l.\label{4.4}
		\end{equation}
	In order to save memory, the initial matrix $H_{k-l}$ is replaced by
		\begin{equation}
	H_{k}^{(0)}=\gamma_{k}{I} . \label{4.5}
		\end{equation}
	Here $\gamma_{k}>0 $ is usually determined by the Barzilai-Borwein method \cite{3} as follows
	\begin{equation}
	\gamma_{k}^{\mathrm{BB} 1}=\frac{\mathbf{y}_{k}^{\top} \mathbf{s}_{k}}{\left\|\mathbf{y}_{k}\right\|^{2}} \quad \text { and } \quad \gamma_{k}^{\mathrm{BB} 2}=\frac{\left\|\mathbf{s}_{k}\right\|^{2}}{\mathbf{y}_{k}^{\top} \mathbf{s}_{k}}.\label{4.6}
	\end{equation}	
	If  $ l\geq k $, $ H_k $ is generated by the traditional BFGS method. The L-BFGS method can be implemented in an inexpensive two-loop recursion way which is shown in Algorithm \ref{alg:1}.
		
	The $\mathbf{p}_{k}=-H_{k}\mathbf{g}\left(\mathbf{x}_{k}\right)$ generated by L-BFGS is a gradient-related descent direction for $H_{k}$ beging a positive definite matrix, where the proof process is presented in Lemma \ref{lemB.2}. Then we find a proper step size along the direction $\mathbf{p}_{k}$ that satisfies the Armijo condition. The computation method for the  model \eqref{3.3} is obtained and we demonstrate it in Algorithm \ref{alg:5}.
	\begin{algorithm}[h]
	\caption{The two-loop recursion for L-BFGS \cite{10,21}.} \label{alg:1}
		\begin{algorithmic}[1]		
			\STATE $\mathbf{q} \leftarrow-\mathbf{g}\left(\mathbf{x}_{k}\right)$
			\FOR { $i=k-1, k-2, \ldots, k-l $ }
			\STATE  $\quad \alpha_{i} \leftarrow \rho_{i} \mathbf{s}_{i}^{\top} \mathbf{q} $
			\STATE  $ \quad \mathbf{q} \leftarrow \mathbf{q}-\alpha_{i} \mathbf{y}_{i} $
			\ENDFOR
			\FOR{for $ i=k-l, k-l+1, \ldots, k-1 $}
			\STATE  $ \quad \beta \leftarrow \rho_{i} \mathbf{y}_{i}^{\top} \mathbf{p} $
			\STATE  $ \quad \mathbf{p} \leftarrow \mathbf{p}+\mathbf{s}_{i}\left(\alpha_{i}-\beta\right) $
			\ENDFOR
			\STATE Stop with result  $\mathbf{p}=-H_{k} \mathbf{g}\left(\mathbf{x}_{k}\right)$.
		\end{algorithmic}
	\end{algorithm}
	\begin{algorithm}[h]
\caption{Low rank triple decomposition for obtaining super-resolution image (TTDSR).}\label{alg:5}
		\begin{algorithmic}[1]
			\STATE Choose constant $r$ and an initial iterate $\mathbf{x}_{0}\in \mathbb{Re}^{m_{1} r r + m_{2} r r + n_{3} r r}$. Select parameters ${l}>0$, $\beta \in \left(0,1\right)$, and $\sigma\in\left(0,1\right)$. Compute $f_{0}=f\left(\mathbf{x}_{0}\right)$ and $\mathbf{g}_{0}=\nabla f\left(\mathbf{x}_{0}\right)$. Set $k \leftarrow 0.$
			\WHILE {the sequence of iterates does not converge}
			\STATE Generate  $\mathbf{p}_{k}=-H_{k} \mathbf{g}\left(\mathbf{x}_{k}\right) $ by Algorithm \ref{alg:1}.
			\STATE Choose the smallest nonnegative integer $\omega$ such that the step size $\alpha=\beta^{\omega }$ satisfies
			\begin{equation}
			f\left(\mathbf{x}_{k}+\alpha \mathbf{p}_{k}\right) \leq f\left(\mathbf{x}_{k}\right)+\sigma \alpha \mathbf{p}_{k}^{\top} \mathbf{g}_{k}.\label{3.30}
			\end{equation}\\
			\STATE Let $\alpha_{k}=\beta^{\omega }$ and update the new iterate $\mathbf{x}_{k+1}=\mathbf{x}_{k}+\alpha_{k} \mathbf{p}_{k}$.
			\STATE Compute $f_{k+1}=f\left(\mathbf{x}_{k+1}\right)$ and $\mathbf{g}_{k+1}=\nabla f\left(\mathbf{x}_{k+1}\right)$.
			\STATE Compute $\mathbf{y}_{k}$, $\mathbf{s}_{k}$ and $\mathbf{\rho}_{k}$ by (\ref{4.1}) and (\ref{4.2}), respectively.
			\STATE $k \leftarrow k+1.$
			\ENDWHILE
		\end{algorithmic}
	\end{algorithm}
\subsection{Gradient calculation}
	At each iteration in Algorithm \ref{alg:5}, we need to compute the gradient of the objective function. The gradient can be calculated via reshaping the tensor variable into vector or matrix form. For small scale problems, the gradient is easier to obtain when the  the tensor variable is turned into a vector than into a matrix, while for large scale problems, the  gradient in vector form always leads to insufficient memory in the computation process. In this subsection, we derive the matrix form of the gradient and the vector form  $\nabla f$ is given in the Appendix \ref{AppendixA}.

	Denote $\mathcal{H}=\mathcal{A} \times_{1} P_{1}\in \mathbb{Re}^{n_{1}\times m \times n}$, $\mathcal{K}=\mathcal{B} \times_{2} P_{2} \in  \mathbb{Re}^{l\times n_{2} \times n}$, and $\mathcal{G}=\mathcal{C} \times_{3} P_{3} \in \mathbb{Re}^{l\times m\times m_{3}}.$  According to (\ref{3.10}) and (\ref{3.2}), we get the following  equations	
	\begin{align}
	\left(Y_{h}\right)_{n_{1}\times{n_{2}{n_{3}}}}&=P_{n_{1} \times m_{1}} A_{m_{1} \times nm} \left(E_{m} \otimes K_{n \times n_{2}l}\right)\left(C_{lm \times n_{3}} \otimes E_{n_{2}}\right),  \label{a.1}\\
	\left(Y_{h}\right)_{n_{2}\times{{n_{3}n_{1}}}}&=P_{n_{2} \times m_{2}} B_{m_{2} \times ln} \left(E_{n} \otimes C_{l \times n_{3}m}\right)\left(H_{mn  \times n_{1}} \otimes E_{n_{3}}\right), \label{a.2}\\
	\left(Y_{h}\right)_{n_{3}\times{{n_{2}n_{1}}}}&=C_{n_{3} \times lm}\left(E_{m} \otimes K_{l \times n_{2}n}\right)\left(H_{nm \times n_{1}} \otimes E_{n_{2}}\right), \label{a.3}\\
	\left(Y_{m}\right)_{m_{1}\times{{m_{2}m_{3}}}}&=A_{m_{1} \times nm}\left(E_{m} \otimes B_{n \times m_{2}l}\right)\left(G_{lm \times m_{3}} \otimes E_{m_{2}}\right),  \label{a.4}\\
	\left(Y_{m}\right)_{m_{2}\times{{m_{3}m_{1}}}}&=B_{m_{2} \times ln}\left(E_{n} \otimes G_{l \times m_{3}m}\right)\left(A_{mn \times m_{1}} \otimes E_{m_{3}}\right),  \label{a.5} \\
	\left(Y_{m}\right)_{m_{3}\times{{ m_{2} m_{1} }}}&=P_{m_{3}\times n_{3}}C_{n_{3}\times lm}\left(E_{m} \otimes B_{l \times m_{2} n}\right)\left(A_{nm \times m_{1}} \otimes E_{m_{2}}\right).  \label{a.6}
	\end{align}
	Define ${P}_{2} {B}_{{m}_{2} \times {ln}}={K}_{{n}_{2} \times ln}$ and $ {P}_{3} {C}_{{n}_{3} \times lm}={G}_{{m}_{3} \times {lm}}.$ The function $f(\mathbf{x})$ in (\ref{3.3}) is transformed into the objective function of matrix variables ${A_{m_{1}\times nm}}$, ${B_{m_{2}\times ln}}$ and ${C_{n_{3}\times lm}}$ as follows
	\begin{equation}\label{3.24}
	\begin{split}
	f(A,B,C)&=\left\|P_{1}A_{m_{1} \times n m}\left(E_{m} \otimes K_{n \times n_{2} l}\right)\left(C_{n_{3} \times l m}^{T} \otimes E_{n_{2}}\right)-\left(Y_{h}\right)_{n_{1} \times n_{2} n_{3}}\right\|_{F}^{2}\\&+\left\|A_{m_{1} \times n m}\left(E_{m} \otimes B_{n \times m_{2} L}\right)\left(\left(P_{3} C_{n_{3} \times lm }\right)^{T} \otimes E_{m_{2}}\right)-\left(Y_{m}\right)_{m_{1} \times m_{2} m_{3}}\right\|_{F}^{2} \\&+\mu\left(\left\|{A}_{{m}_{1} \times {nm}}\right\|_{{F}}^{2}+\left\|{B}_{{m}_{2} \times {ln}}\right\|_{{F}}^{2}+\left\|{C}_{{n}_{3} \times {lm}}\right\|_{{F}}^{2}\right),
	\end{split}
	\end{equation}
		\begin{equation}\label{3.25}
		\begin{split}
		f(A,B,C)
		&=\left\|P_{2}B_{m_{2} \times ln }\left(E_{n} \otimes C_{l \times n_{3} m}\right)\left(\left(P_{1}{A_{m_{1}\times mn}}\right)^{T} \otimes E_{n_{3}}\right)-\left(Y_{h}\right)_{n_{2} \times n_{3} n_{1}}\right\|_{F}^{2}\\
		&+\left\|B_{m_{2} \times ln}\left(E_{n} \otimes G_{l \times m_{3} m}\right)\left(\left(A_{m_{1}\times mn} \right)^{T} \otimes E_{m_{3}}\right)-\left(Y_{m}\right)_{m_{2} \times m_{3} m_{1}}\right\|_{F}^{2} \\
		&+\mu\left(\left\|{A}_{{m}_{1} \times {nm}}\right\|_{{F}}^{2}+\left\|{B}_{{m}_{2} \times{ln}}\right\|_{{F}}^{2}+\left\|{C}_{{n}_{3} \times {lm}}\right\|_{{F}}^{2}\right),
		\end{split}
		\end{equation}
	\begin{equation}\label{3.26}
	\begin{split}
	f(A,B,C)
	&=\left\|C_{n_{3}\times lm}\left(E_{m} \otimes K_{l \times n_{2} n}\right)\left(\left(P_{1}{A_{m_{1}\times mn}}\right)^{T} \otimes E_{n_{2}}\right)-\left(Y_{h}\right)_{n_{3} \times n_{2} n_{1}}\right\|_{F}^{2}\\
	&+\left\|P_{3}C_{n_{3} \times lm}\left(E_{m} \otimes B_{l \times m_{2} n}\right)\left(\left(A_{m_{1}\times nm} \right)^{T} \otimes E_{m_{2}}\right)-\left(Y_{m}\right)_{m_{3} \times m_{2} m_{1}}\right\|_{F}^{2}\\
	&+\mu\left(\left\|{A}_{{m}_{1} \times {nm}}\right\|_{F}^{2}+\left\|{B}_{{m}_{2} \times {ln}}\right\|_{{F}}^{2}+\left\|{C}_{{n}_{3} \times {l m}}\right\|_{{F}}^{2}\right).
	\end{split}
	\end{equation}
	For any matrix $X,$ we have $\|X\|_F^2=tr(X^{\top}X).$ Therefore we review the derivatives of some useful trace functions with respect to $X:$
	\begin{equation}\label{3.28}
	\normalsize  \begin{split}
	&\dfrac{\partial tr(AXBX^{T}C)}{\partial X}=A^{T}C^{T}XB^{T}+CAXB,\quad\quad
	\dfrac{\partial tr(AX^{T}B)}{\partial X}=BA,\quad\\
	&\dfrac{\partial tr(ABA^{T})}{\partial A}=A(B+B^{T}),\quad\quad\quad\quad\quad\quad\quad\quad\quad
	\dfrac{\partial tr(BA^{T})}{\partial A}=B,\quad\\
	&\dfrac{\partial tr(AXB)}{\partial X}=A^{T}B^{T},\qquad\qquad\qquad\qquad\qquad\quad\quad
	\dfrac{\partial tr(AB)}{\partial A}=B^{T}.\quad
	\end{split}
	\end{equation}
	For simplicity, we denote
	\begin{align*}
{E_{1}}&=\left(E_{m} \otimes B_{n \times m_{2} l}\right)\left(\left(P_{3} C_{n_{3} \times lm }\right)^{T} \otimes E_{m_{2}}\right),
&{D_{1}}&=\left(E_{m} \otimes K_{n \times n_{2} l}\right)\left(C_{n_{3} \times lm}^{T} \otimes E_{n_{2}}\right),  \\
	{F_{1}}&=\left(E_{n} \otimes C_{l \times n_{3} m}\right)\left(\left(P_{1}{A_{m_{1}\times mn}}\right)^{T} \otimes E_{n_{3}}\right),  &{G_{1}}&=\left(E_{n} \otimes G_{l \times m_{3} m}\right)\left(A_{m_{1}\times mn}^{T} \otimes E_{m_{3}}\right),\\
	H_{1}&=\left(E_{m}\otimes  K_{l\times n_{2} n}\right)\left(\left(P_{1}{A_{m_{1}\times mn}}\right)^{T} \otimes E_{n_{2}}\right),&K_{1}&=\left(E_{m} \otimes B_{l \times m_{2} n}\right)\left(A_{m_{1}\times nm}^{T} \otimes E_{m_{2}}\right).
	\end{align*}
	Thus,  the partial derivatives with respect to $A,$ $B,$ and $C$ are
		\begin{equation}\label{4.31}
		\begin{split}
		\dfrac{\partial f({A},{B},{C})}{\partial {A}}
		=&\dfrac{\partial \left(\left\|{P_{1}}{A}{D_{1}}-{Y_{h}}\right\|_{F}^{2}+ \left\|{A}{E_{1}}-{Y_{m}}\right\|_{F}^{2}+\mu\left(\left\|{A} \right\|_{{F}}^{2}+\left\|{B}\right\|_{{F}}^{2}+\left\|{C}\right\|_{{F}}^{2}\right)\right)}{\partial {A}}\\\vspace{1.5ex}
		=&\dfrac{\partial tr\left( \langle {P_{1}}{A}{D_{1}}-{Y_{h}}, {P_{1}}{A}{D_{1}}-{Y_{h}} \rangle + \langle {A}{E_{1}}-{Y_{m}}, {A}{E_{1}}-{Y_{m}}\rangle  \right)}{\partial {A}}+ 2\mu{A}\\
		=&\dfrac{\partial tr\left({P_{1}}{A}{D_{1}}{D_{1}^{T}}{A^{T}}{P_{1}^{T}}-{P_{1}}{A}{D_{1}}{Y_{h}^{T}}-{Y_{h}}{D_{1}^{T}}{A^{T}}{P_{1}^{T}}-{Y_{h}}{Y_{h}^{T}}\right)}{\partial {A}}\\
		+&\dfrac{\partial tr \left({A}{E_{1}}{E_{1}^{T}}{A^{T}}-{A}{E_{1}}{Y_{m}^{T}}-{Y_{m}}{E_{1}^{T}}{A^{T}}+{Y_{m}^{T}} {Y_{m}}\right)}{{\partial {A}}}+2\mu {A}\\
		=&2\left({P_{1}^{T}}{P_{1}}A{D_{1}}{D_{1}^{T}}-{P_{1}^{T}}{Y_{h}}{D_{1}^{T}}+{A}{E_{1}}{E_{1}^{T}}-{Y_{m}}{E_{1}^{T}}+\mu {A}\right),
		\end{split}
		\end{equation}
		\begin{equation}\label{4.32}
		\begin{split}
		\dfrac{\partial f({A},{B},{C})}{\partial {B}}
		=&\dfrac{\partial \left(\left\|{P_{2}}{B}{F_{1}}-{Y_{h}}\right\|_{F}^{2}+ \left\|{B}{G_{1}}-{Y_{m}}\right\|_{F}^{2}+\mu\left(\left\|{A} \right\|_{{F}}^{2}+\left\|{B}\right\|_{{F}}^{2}+\left\|{C}\right\|_{{F}}^{2}\right)\right)}{\partial {B}}\\
		=&\dfrac{\partial tr\left( \langle {P_{2}}{B}{F_{1}}-{Y_{h}}, {P_{2}}{B}{F_{1}}-{Y_{h}}\rangle + \langle {B}{G_{1}}-{Y_{m}}, {B}{G_{1}}-{Y_{m}}\rangle  \right)}{\partial {B}}+ 2\mu{B}\\
		=&\dfrac{\partial tr\left({P_{2}}{B}{F_{1}}{F_{1}^{T}}{B^{T}}{P_{2}^{T}}-{P_{2}}{B}{F_{1}}{Y_{h}^{T}}-{Y_{h}}{F_{1}^{T}}{B^{T}}{P_{2}^{T}}-{Y_{h}}{Y_{h}^{T}}\right)}{\partial {B}}\\
		+&\dfrac{\partial tr \left({B}{G_{1}}{G_{1}^{T}}{B^{T}}-{B}{G_{1}}{Y_{m}^{T}}-{Y_{m}}{G_{1}^{T}}{B^{T}}+{Y_{m}^{T}} {Y_{m}}\right)}{{\partial {B}}}+2\mu {B}\\
		=&2\left({P_{2}^{T}}{P_{2}}B{F_{1}}{F_{1}^{T}}-{P_{2}^{T}}{Y_{h}}{F_{1}^{T}}+{B}{G_{1}}{G_{1}^{T}}-{Y_{m}}{G_{1}^{T}}+\mu {B}\right),
		\end{split}
		\end{equation}
		\begin{equation}\label{4.33}
		\begin{split}
		\dfrac{\partial f({A},{B},{C})}{\partial {C}}
		=&\dfrac{\partial \left( \left\|{C}{H_{1}}-{Y_{h}}\right\|_{F}^{2}+\left\|{P_{3}}{C}{K_{1}}-{Y_{m}}\right\|_{F}^{2}+\mu\left(\left\|{A} \right\|_{{F}}^{2}+\left\|{B}\right\|_{{F}}^{2}+\left\|{C}\right\|_{{F}}^{2}\right)\right)}{\partial {C}}\\
		=&\dfrac{\partial tr \left(\langle  {C}{H_{1}}-{Y_{h}}, {C}{H_{1}}-{Y_{h}}\rangle  + \langle {P_{3}}{C}{K_{1}}-{Y_{m}}, P_{3}{C}{K_{1}}-{Y_{m}}\rangle  \right)}{\partial {C}}+ 2\mu{C}\\
		=&\dfrac{\partial tr\left({C}{H_{1}}{H_{1}^{T}}{C^{T}}-{C}{H_{1}}{Y_{h}^{T}}-{Y_{h}}{H_{1}^{T}}{C^{T}}+{Y_{h}^{T}} {Y_{h}}\right)}{\partial {C}}\\
		+&\dfrac{\partial tr\left({P_{3}}{C}{K_{1}}{K_{1}^{T}}{C^{T}}{P_{3}^{T}}-{P_{3}}{C}{K_{1}}{Y_{m}^{T}}-{Y_{m}}{K_{1}^{T}}{C^{T}}{P_{3}^{T}}-{Y_{m}}{Y_{m}^{T}}\right)}{\partial {C}}+2\mu {C}\\
		=&2\left({C}{H_{1}}{H_{1}^{T}}-{Y_{h}}{H_{1}^{T}}+{P_{3}^{T}}{P_{3}}C{K_{1}}{K_{1}^{T}}-{P_{3}^{T}}{Y_{m}}{K_{1}^{T}}+\mu {C}\right).
		\end{split}
		\end{equation}
\section{Convergence analysis}
\label{section:5}
In this section, we analyze the convergence of $ \|\mathbf{g}(\mathbf{x})\|$ and show that our proposed method produces a globally convergent iteration.
\begin{lemma}	\label{lem5.1}
	There exists a positive number $ M $ such that
	\begin{equation}
	\lvert f(\mathbf{x}) \rvert\leq M,  \|\mathbf{g}(\mathbf{x})\| \leq M,  \|H(\mathbf{x})\| \leq M.
	\end{equation}
\end{lemma}
\begin{proof}
	Because $\mathbf{p}_{k}=-H_{k} \mathbf{g}\left(\mathbf{x}_{k}\right) $ generated by L-BFGS is a gradient-related descent direction. There exists a positive number satisfies the Armijo condition in (\ref{3.30}) and the sequence of objective function  value $ \left\{f\left(\mathbf{x}_{k}\right)\right\}  $  decreases monotonically. Since the value of the objective function  $f\left(\mathbf{x}\right)$ in (\ref{3.3}) is nonnegative, $ f\left(\mathbf{x}\right)$ is bounded, i.e. $0 \leq f(\mathbf{x}) \leq M$ holds.
	
	The regularization term   $\mu\left(\|\mathcal{A}\|_{F}^{2}+\|\mathcal{B}\|_{F}^{2}+\|\mathcal{C}\|_{F}^{2}\right)$ in the objective function indicates that the sequence $\left\{\mathbf{x}_{k}\right\}$ is bounded.
	Because $\mathbf{g}\left(\mathbf{x}\right)$ is a linear function of  $\mathbf{x}$, $\|\mathbf{g}(\mathbf{x})\|$ is also bounded. The proof of the boundedness of $ H(\mathbf{x}) $ is shown  in Appendix \ref{AppendixB}. Thus, we get Lemma \ref{lem5.1}.
\end{proof}

Furthermore,  we have the following Lemma \ref{lem5.2}.
\begin{lemma}	\label{lem5.2}
	There exist constants $C_{m}$ and $C_{M}$ satisfy
	$ 0<C_{m} \leq 1 \leq C_{M} $,
	\begin{equation}\label{5.2}
	\mathbf{p}_{k}^{\top} \mathbf{g}\left(\mathbf{x}_{k}\right) \leq-C_{m}\left\|\mathbf{g}\left(\mathbf{x}_{k}\right)\right\|^{2} \quad \text { and } \quad\left\|\mathbf{p}_{k}\right\| \leq C_{M}\left\|\mathbf{g}\left(\mathbf{x}_{k}\right)\right\|.
	\end{equation}
\end{lemma}
\begin{proof}
	The proof can be found in the Appendix \ref{AppendixB}.
\end{proof}
Because of boundedness and monotonicity of $\{f\left(\mathbf{x}_{k}\right) \}$, the sequence of function value converges. The conclusion is given in Theorem \ref{thm5.4}.
\begin{theorem}\label{thm5.4}
	Assume that Algorithm $\ref{alg:5}$ generates an infinite sequence of function values $ \left\{f\left(\mathbf{x}_{k}\right)\right\}$. Then, there exists a constant ${f}_{*}$ such that
	\begin{equation}
	\lim _{k\rightarrow \infty} f\left(\mathbf{x}_{k}\right)=f_{*}.
	\end{equation}
\end{theorem}

Next, we prove that every accumulation point of iterates $ \{\mathbf{x}_{k}\} $ is a first-order stationary point. At last, by utilizing the Kurdyka-$\L$ojasiewicz property \cite{2}, we show that the sequence of iterates $\left\{\mathbf{x}_{k}\right\}$ is also convergent. The following lemma means that the step size is lower bounded.

\begin{lemma}	\label{lem5.5}
	There exists a constant $ \alpha_{\min }>0 $ such that
	\begin{equation}
	\alpha_{k} \geq \alpha_{\min }>0,  \quad \quad \forall k\in \mathbf{N}^{+}.
	\end{equation}
\end{lemma}
\begin{proof}
	Let $0<\alpha \leq \tilde{\alpha}=\frac{(1-\sigma) C_{m}}{\frac{1}{2} M C_{M}^{2}}$.
	From Lemma \ref{lem5.2}, for $\alpha \in (0,\tilde{\alpha}]$, it yields that
	\begin{equation}\label{5.5}
	\begin{split}
	\alpha \mathbf{p}_{k}^{T} \mathbf{g}_{k}+\frac{1}{2} M \alpha^{2}\left\|\mathbf{p}_{k}\right\|^{2}-\sigma\alpha \mathbf{p}_{k}^{T} \mathbf{g}_{k}
	&=(1-\sigma)\alpha\mathbf{p}_{k}^{T} \mathbf{g}_{k}+\frac{1}{2} M \alpha^{2}\left\|\mathbf{p}_{k}\right\|^{2}\\
	&\leq (1-\sigma)\alpha(-C_{m}\left\|\mathbf{g}\left(\mathbf{x}_{k}\right)\right\|^{2})+\frac{1}{2} M \alpha^{2}C_{M}^{2}\left\|\mathbf{g}_{k}\right\|^{2}\\
	&\leq 0.
	\end{split}
	\end{equation}
	From Taylor's mean value theorem  and $ \mathbf{x}_{k+1}(\alpha)=\mathbf{x}_{k}+\alpha \mathbf{p}_{k}$, we have
	\begin{equation}\label{4.8}
	\begin{aligned}
	f\left(\mathbf{x}_{k+1}(\alpha)\right)-f\left(\mathbf{x}_{k}\right) &\leq \mathbf{g}_{k}^{T}\left(\mathbf{x}_{k+1}(\alpha)-\mathbf{x}_{k}\right)+\frac{1}{2} M\left\|\mathbf{x}_{k+1}(\alpha)-\mathbf{x}_{k}\right\|^{2} \\
	&=\alpha \mathbf{p}_{k}^{T} \mathbf{g}_{k}+\frac{1}{2} M \alpha^{2}\left\|\mathbf{p}_{k}\right\|^{2}\\
	&\leq \sigma\alpha \mathbf{p}_{k}^{T} \mathbf{g}_{k},
	\end{aligned}
	\end{equation}
	where the last inequality is valid owing to (\ref{5.5}). The rule of  inexact line search indicates $\alpha_{k} \geq \alpha_{\min }=\tilde{\alpha}\beta.$ Hence, we find out a lower bound on the step size $\alpha$.
\end{proof}

The following theorem  proves that every accumulation point of iterates  $\left\{\mathbf{x}_{k}\right\}$ is a first-order stationary point.
\begin{theorem} \label{thm5.5}
	Suppose that  Algorithm $\ref{alg:5}$ generates an infinite sequence of iterates  $\left\{\mathbf{x}_{k}\right\} $. Then,
	\begin{equation}
	\lim _{k \rightarrow \infty}\left\|\mathbf{g}\left(\mathbf{x}_{k}\right)\right\|=0.
	\end{equation}
\end{theorem}
\begin{proof}
	From Lemma \ref{lem5.2} and (\ref{3.30}), we get
	\begin{equation}\label{eq5.8}
	f\left(\mathbf{x}_{k}\right)-f\left(\mathbf{x}_{k+1}\right) \geq-\sigma \alpha_{k} \mathbf{p}_{k}^{\top} \mathbf{g}\left(\mathbf{x}_{k}\right) \geq \sigma \alpha_{k} C_{m}\left\|\mathbf{g}\left(\mathbf{x}_{k}\right)\right\|^{2}.
	\end{equation}
	The functional series $\sum_{k=1}^{\infty}\left[f\left(\mathbf{x}_{k}\right)-f\left(\mathbf{x}_{k+1}\right)\right]$ satisfies
	\begin{equation}
	\begin{split}
	2 M \geq f\left(\mathbf{x}_{1}\right)-f_{*}&=\sum_{k=1}^{\infty}\left[f\left(\mathbf{x}_{k}\right)-f\left(\mathbf{x}_{k+1}\right)\right]\\ &\geq \sum_{k=1}^{\infty} \sigma \alpha_{k} C_{m}\left\|\mathbf{g}\left(\mathbf{x}_{k}\right)\right\|^{2} \\
&\geq \sum_{k=1}^{\infty} \sigma \alpha_{\min } C_{m}\left\|\mathbf{g}\left(\mathbf{x}_{k}\right)\right\|^{2}.
	\end{split}
	\end{equation}
	That is to say,
	\begin{equation}
	\sum_{k=1}^{\infty}\left\|\mathbf{g}\left(\mathbf{x}_{k}\right)\right\|^{2} \leq \frac{2 M}{\sigma \alpha_{\min } C_{m}}<+\infty.
	\end{equation}
	Hence, $\left\|\mathbf{g}\left(\mathbf{x}_{k}\right)\right\|$ converges to zero.
\end{proof}

Analysis of proximal methods for nonconvex and nonsmooth optimization frequently uses the Kurdyka- $\L$ojasiewicz (KL) property \cite{2}. Since the objective function $f(\mathbf{x}) $ is a polynomial and the KL property below holds, we use the KL property to prove the convergence of algorithm.
\begin{proposition}{(Kurdyka-$\L$ojasiewicz (KL) property)}
	Suppose that $\mathbf{x}_{*} $ is a stationary point of $ f(\mathbf{x})$.  There is a neighborhood
	$ \mathscr{U} $ of  $ \mathbf{x}_{*} $, an exponent  $\theta \in[0,1) $ and a positive constant $ K $ such that for all $\mathbf{x}\in \mathscr{U}$, the following inequality holds:
	\begin{equation}
	\left|f(\mathbf{x})-f\left(\mathbf{x}_{*}\right)\right|^{\theta} \leq K \|\mathbf{g}(\mathbf{x})\|.
	\end{equation}
	In particular, we define $0^{0}\equiv 0 $.
\end{proposition}
\begin{lemma} \label{lem5.7}
	Suppose that  $\mathbf{x}_{*} $ is a stationary point of  $f(\mathbf{x})$ and
	$\mathscr{A}\left(\mathbf{x}_{*},\delta\right)=\{\mathbf{x} \in \left.\mathbb{R}^{n}:\left\|\mathbf{x}-\mathbf{x}_{*}\right\| \leq \delta\right\} \subseteq \mathscr{U} $ is a neighborhood of $\mathbf{x}_{*}.$
	Let $\mathbf{x}_{1} $ be an initial point satisfying
	\begin{equation}
	\delta>\frac{ C_{M} K}{\sigma C_{m}(1-\theta)}\left|f\left(\mathbf{x}_{1}\right)-f\left(\mathbf{x}_{*}\right)\right|^{1-\theta}+\left\|\mathbf{x}_{1}-\mathbf{x}_{*}\right\|.
	\end{equation}
	Then, the following assertions hold:
	\begin{equation}\label{4.15}
	\mathbf{x}_{k} \in \mathscr{A}\left(\mathbf{x}_{*}, \delta\right), \quad k=1,2, \ldots
	\end{equation}
	and
	\begin{equation}\label{4.16}
	\sum_{k=1}^{\infty}\left\|\mathbf{x}_{k+1}-\mathbf{x}_{k}\right\| \leq \frac{C_{M} K}{\sigma C_{m}(1-\theta)}\left|f\left(\mathbf{x}_{1}\right)-f\left(\mathbf{x}_{*}\right)\right|^{1-\theta}.
	\end{equation}
\end{lemma}
\begin{proof}
	The theorem is proved by induction. Obviously, $ \mathbf{x}_{1} \in \mathscr{A}\left(\mathbf{x}_{*}, \delta\right) $.
	Now, we assume $\mathbf{x}_{i} \in \mathscr{A}\left(\mathbf{x}_{*}, \delta\right) $ for all  $i=1, \ldots, k $ and KL property holds at these points. Define a concave function
	\begin{equation}
	\phi(q) \equiv \frac{K}{1-\theta}\left|q-f\left(\mathbf{x}_{*}\right)\right|^{1-\theta}\quad and \quad q>f\left(\mathbf{x}_{*}\right).
	\end{equation}
	Its derivative function is
	\begin{equation}
	\phi^{\prime}(q)=\frac{K}{\left|q-f\left(\mathbf{x}_{*}\right)\right|^{\theta}}.
	\end{equation}
	For $ i=1, \ldots, k $,   the first-order condition of the concave function $\phi^{\prime}(q)$ at $f(\mathbf{x}_{i})$ is
	\begin{equation}\label{5.19}
	\begin{aligned}
	\phi\left(f\left(\mathbf{x}_{i}\right)\right)-\phi\left(f\left(\mathbf{x}_{i+1}\right)\right) & \geq \phi^{\prime}\left(f\left(\mathbf{x}_{i}\right)\right)\left(f\left(\mathbf{x}_{i}\right)-f\left(\mathbf{x}_{i+1}\right)\right) \\
	&=\frac{K}{\left|f\left(\mathbf{x}_{i}\right)-f\left(\mathbf{x}_{*}\right)\right|^{\theta}}\left(f\left(\mathbf{x}_{i}\right)-f\left(\mathbf{x}_{i+1}\right)\right).
	\end{aligned}
	\end{equation}
	The equation  (\ref{5.19}) and KL property mean
	\begin{equation} \phi\left(f\left(\mathbf{x}_{i}\right)\right)-\phi\left(f\left(\mathbf{x}_{i+1}\right)\right)\geq\frac{1}{\left\|\mathbf{g}\left(\mathbf{x}_{i}\right)\right\|}\left(f\left(\mathbf{x}_{i}\right)-f\left(\mathbf{x}_{i+1}\right)\right).
	\end{equation}
	By Lemma \ref{lem5.2} and (\ref{eq5.8}), we have
	\begin{equation}\label{5.21}
	\phi\left(f\left(\mathbf{x}_{i}\right)\right)-\phi\left(f\left(\mathbf{x}_{i+1}\right)\right)\geq \sigma \alpha_{i} C_{m}\left\|\mathbf{g}\left(\mathbf{x}_{i}\right)\right\| \geq \frac{\sigma C_{m}}{ C_{M}}\left\|\mathbf{x}_{i+1}-\mathbf{x}_{i}\right\|,
	\end{equation}
	where the last inequality is valid because
	\begin{equation}
	\left\|\mathbf{x}_{k+1}-\mathbf{x}_{k}\right\| = \alpha_{k}\left\|\mathbf{p}_{k}\right\| \leq \alpha_{k} C_{M}\left\|\mathbf{g}\left(\mathbf{x}_{k}\right)\right\|.
	\end{equation}
	The upper bound of $\|\mathbf{x}_{k+1}-\mathbf{x}_{*}\|$ is
	\begin{equation}\label{4.23}
	\begin{aligned}
	\left\|\mathbf{x}_{k+1}-\mathbf{x}_{*}\right\| & \leq \sum_{i=1}^{k}\left\|\mathbf{x}_{i+1}-\mathbf{x}_{i}\right\|+\left\|\mathbf{x}_{1}-\mathbf{x}_{*}\right\| \\
	& \leq \frac{C_{M}}{\sigma C_{m}} \sum_{i=1}^{k} \phi\left(f\left(\mathbf{x}_{i}\right)\right)-\phi\left(f\left(\mathbf{x}_{i+1}\right)\right)+\left\|\mathbf{x}_{1}-\mathbf{x}_{*}\right\| \\
	& \leq \frac{C_{M}}{\sigma C_{m}} \left(\phi\left(f\left(\mathbf{x}_{1}\right)\right)-\phi\left(f\left(\mathbf{x}_{k+1}\right)\right)\right)+\left\|\mathbf{x}_{1}-\mathbf{x}_{*}\right\| \\
	& \leq \frac{C_{M}}{\sigma C_{m}} \phi\left(f\left(\mathbf{x}_{1}\right)\right)+\left\|\mathbf{x}_{1}-\mathbf{x}_{*}\right\|\\
	&<\delta,
	\end{aligned}
	\end{equation}
	which means $ \mathbf{x}_{k+1} \in \mathscr{A}\left(\mathbf{x}_{*}, \delta\right)$  and (\ref{4.15}) holds. Moreover, according to (\ref{5.21}), we obtain
	\begin{equation}\label{5.24}
	\sum_{k=1}^{\infty}\left\|\mathbf{x}_{k+1}-\mathbf{x}_{k}\right\| \leq \frac{C_{M}}{\sigma C_{m}} \sum_{k=1}^{\infty} \phi\left(f\left(\mathbf{x}_{k}\right)\right)-\phi\left(f\left(\mathbf{x}_{k+1}\right)\right) \leq \frac{C_{M}}{\sigma C_{m}} \phi\left(f\left(\mathbf{x}_{1}\right)\right).
	\end{equation}
	Thus, the proof of (\ref{4.16}) is complete.
\end{proof}
The sequence of iterates  $\left\{\mathbf{x}_{k}\right\} $ is demonstrated to converge to a unique accumulation point next.
\begin{theorem}
	Suppose that Algorithm $\ref{alg:5}$  generates an infinite sequence of iterates  $\left\{\mathbf{x}_{k}\right\} $. $\left\{\mathbf{x}_{k}\right\} $ converges to a unique first-order stationary point $\mathbf{x}_{*} $, i.e.
	\begin{equation}\label{5.25}
	\lim _{k \rightarrow \infty} \mathbf{x}_{k}=\mathbf{x}_{*}.
	\end{equation}
\end{theorem}

\begin{proof}
	Clearly, according to (\ref{5.24}), the sequence $\{\mathbf{x}_{k}\}$ satisfies
	\begin{equation}\label{5.26}
	\sum_{k=1}^{\infty}\left\|\mathbf{x}_{k+1}-\mathbf{x}_{k}\right\|<+\infty
	\end{equation}
	and is a Cauchy sequence.
	Owing to the boundedness of $\{\mathbf{x}_{k}\}$, there exists an accumulate point $\mathbf{x}_{*}$ of iterates  $\{\mathbf{x}_{k}\}$.  Thus (\ref{5.25}) holds.
\end{proof}

\section{Numerical experiments} 	
\label{section:6}
In this section we demonstrate the performance of our proposed TTDSR method on two datasets. The method is implemented with parameters $l=5$, $\sigma=0.01$, $\beta=0.5$, $\mu=1$. The stopping criteria is
$$	\left\|\mathbf{g}\left(\mathbf{x}_{k}\right)\right\|_{\infty}<10^{-10} $$
or
$$	\left\|\mathbf{x}_{k+1}-\mathbf{x}_{k}\right\|_{\infty}<10^{-16} \quad \text { and } \quad {\left|f\left(\mathbf{x}_{k+1}\right)-f\left(\mathbf{x}_{k}\right)\right|}<10^{-2}.$$
The maximum number of iteration is set to 400. All simulations are run on a HP notebook with 2.5 GHz Intel Core i5 and 4 GB RAM. We use tensorlab 3.0 \cite{26} for basic tensor operations.

In the numerical experiments,  the groundtruth image SRI is artificially degraded to HSI and MSI based on  $ {P}_{1}$, ${P}_{2} $ and $ {P}_{3}.$  For matrices $P_1$ and $P_2$ of spatial degradation, we follow the Wald's protocol \cite{36} and the degradation process from SRI to HSI is a combination of spatial blurring Gaussian kernel and downsampling. The downsampling and Gaussian kernel have parameters $d$ and $q$, respectively. It is common to set the downsampling ratio  $d=4$   and  $q=9$ in the Gaussian kernel. In the following experiments, we also conduct simulations under situations such as $d=6$ and $q=9$, $d=4$ and $q=5$ respectively. In order to obtain MSI from SRI, we  generate the spectral degradation matrix $P_{3}$ through spectral specifications, which are taken from LANDSAT or QuickBird specifications of multispectral sensors. The Indian pines and Salinas-A scene datasets are available online at \cite{14}. For Indian pines, the groundtruth image SRI is degraded with the former sensor, while the SRI of Salinas-A scene is degraded with the latter as in \cite{28}. The dimensions of HSI, MSI, SRI images are demonstrated in Table \ref{tab:1}.
\begin{table}
	\begin{tabular}{cccc}
	\toprule
	\text{Image name} &\text {SRI}  & \text { HSI } & \text {MSI}\\
	\midrule
	\text{Indian pines, d=4} & $144\times 144\times 200 $ & $36\times 36\times 200$ & $144\times 144\times 6$ \\
    \text{Indian pines, d=6} & $144\times 144\times 200 $ & $24\times 24\times 200$ & $144\times 144\times 6$ \\
	\text{Salinas-A scene} & $80\times80\times200$ & $20\times20\times200$ & $80\times80\times4$\\
	\bottomrule
	\end{tabular}
\caption{Image size for the SRI experiments.}\label{tab:1}
\end{table}
\subsection{Comparison with other algorithms}
In this subsection we  compare the proposed algorithm with state-of-the-art approaches, including HySure \cite{31} and FUSE \cite{35}, which are based on matrix decompositions. Furthermore, tensor CP \cite{20} and tucker decomposition \cite{28} methods  are also considered. The HySure method is about a convex formulation for SRI via subspace-based regularization proposed by Simo$\acute{e}$s et al, while FUSE describes fast fusion of multiband images based on solving a Sylvester equation proposed by Wei et al. We calculate the following metrics to evaluate the effect of image fusion, which includes re-constructed signal-to-noise ratio (R-SNR), correlation koeffizient (CC), spectral angle mapper (SAM) and the relative dimensional global error (ERGAS) used in \cite{28}. R-SNR and CC are given by
\begin{equation}
\text{R-SNR}=10 \log_{10}\left(\frac{\|\mathcal{Z}\|_{F}^{2}}{\|\widehat{\mathcal{Z}}-\mathcal{Z}\|_{F}^{2}}\right)
\end{equation}
and
\begin{equation}
\text{CC}=\frac{1}{n_3}\left(\sum_{k=1}^{n_3} \rho\left(\mathcal{Z}_{:,:,k}, \hat{\mathcal{Z}}_{:,:, k}\right)\right),
\end{equation}
where $\rho(\cdot, \cdot) $ is the pearson correlation coefficient between the original and  estimated spectral slices.
The metric SAM is
\begin{equation}
\text{SAM}= \frac{1}{m_{1} m_{2}} \sum_{n=1}^{m_{1} m_{2}} \arccos \left(\frac{{Z}_{n,:}^{(3)^{\top}} \widehat{{Z}}_{n,:}^{(3)}}{\left\|{Z}_{n,:}^{(3)}\right\|_{2}\left\|\hat{{Z}}_{n,:}^{(3)}\right\|_{2}}\right)
\end{equation}
and	 calculates the angle between the  original and  estimated spectral fiber.
The performance measurement ERGAS is
\begin{equation}
\text { ERGAS }=\frac{100}{d} \sqrt{\frac{1}{m_{1}m_{2}n_{3}} \sum_{k=1}^{n_{3}} \frac{\left\|\widehat{\mathcal{Z}}_{:,:, k}-\mathcal{Z}_{:,:, k}\right\|_{F}^{2}}{\mu_{k}^{2}}},
\end{equation}
where  $\mu_{k}^{2}$  is the mean value of  $\hat{\mathcal{Z}}_{:,:, k}$. It represents the relative dimensionless global error between SRI and the estimated one. It is the root mean-square error averaged by the size of the SRI.
\begin{table}
	\begin{tabular}{ll|l|ccccc}
		\hline
		& & \multirow{2}*{Algorithm} & \multicolumn{5}{c}{quality evaluation metrics} \\ \cline{4-8}
		& & &\text { R-SNR }  & \text { CC } & \text { SAM } &\text { ERGAS  }& \text { TIME(s) } \\
		\hline
		& &\text { best value } & $+\infty$  & 1 & 0 & 0 & - \\ \hline
        & &\text { STERTO} & 24.8691 & 0.8335 & 2.8220 & 1.2812 & 2.3025 \\
        &d=4&\text { SCOTT } & 16.4046 & 0.7617 & 7.2446 & 2.3651 & 0.5865 \\
		& &\text { HySure } &18.9055 &0.6971 &5.7052 &2.3045 & 39.9998\\
		&q=9&\text { FUSE} &10.3359 &0.6126 &13.8561 &4.5692& 0.3635 \\
		& &\text { TTDSR } & 17.0350 & 0.6712 & 6.7452 & 3.1165 & 2.3624 \\
		\hline
        & &\text { STERTO  } & 23.7512&0.7874&	3.2923	&0.9875&	2.1344\\
        &d=6&\text { SCOTT } & 17.1569&	0.7414&	6.7436&	1.5801&	0.8730\\
		&&\text { HySure } &17.8228&	0.6879&	6.4205&	1.6793&	38.6143\\
		&q=9&\text { FUSE} &11.9157&0.6134&	11.8385&2.8082&	0.3006\\
		&&\text { TTDSR } &17.0350 & 0.6712 & 6.7452 &2.0777&2.1966 \\
		\hline
        &&\text { STERTO  } & 24.8723 & 0.8338 & 2.8213 & 1.2791 & 1.7986 \\
        &d=4&\text { SCOTT } & 15.5634 & 0.7218 & 7.6436 & 2.9810 & 0.5197\\
		&&\text { HySure } &16.8333 &0.6729 &6.9514 &2.8192 & 37.9540\\
		&q=5&\text { FUSE} &9.5533 &0.5678 &14.6419 &5.8727& 0.3169 \\
		&&\text { TTDSR } &17.0350 & 0.6712 & 6.7452 & 3.1165 & 2.386 \\
		\hline
	\end{tabular}
\caption{Comparison of performance of different algorithms on Indian pines.}\label{tab:2}
\end{table}
\begin{figure}[h]
	\centering
	\subfigure[\label{fig:a}]{
		\includegraphics[scale=0.68]{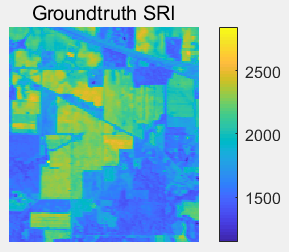}}
	\subfigure[\label{fig:c}]{
		\includegraphics[scale=0.68]{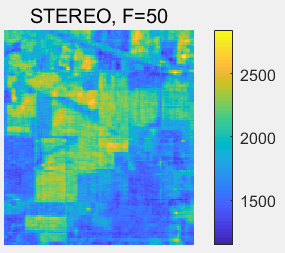}}
	\subfigure[\label{fig:e}]{
		\includegraphics[scale=0.68]{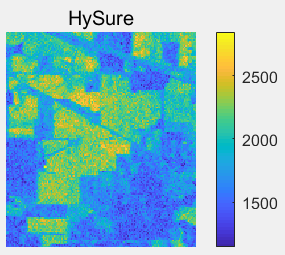}}
	\subfigure[\label{fig:b}]{
		\includegraphics[scale=0.68]{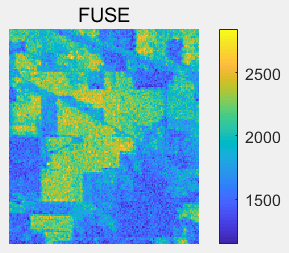}}
\subfigure[\label{fig:b}]{
		\includegraphics[scale=0.68]{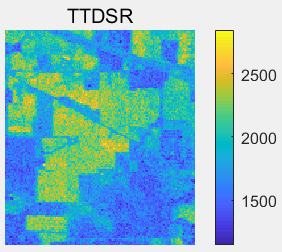}}
	\caption{Spectral slice 120 of the SRI, Indian pines.}
	\label{fig3}
\end{figure}

In practical applications, it is common that the hyperspectral and multispectral images generated by special sensors are noisy. Therefore in the experiments, we  add white Gaussian  noise to HSI and MSI. The first experiment is performed using Indian pines dataset  from hyperspectral remote sensing data platform \cite{14}. White Gaussian noise to $\mathcal{Y}_{h}$ is 21dB, while to $\mathcal{Y}_{m}$ is 25dB. The results are presented in Table \ref{tab:2} and Figure \ref{fig3}. The rank of tucker decomposition in SCOTT is $ \left [70,70, 6 \right]$. According to \cite{20}, tensor rank $F = 50$ of the STEREO method often yields good performance. For HySure method, `E' represents groundtruth number of materials and is set to $16,$ which is chosen as the number of endmembers as \cite{20}. In Table \ref{tab:2}, when $d=4$ and $q=9$, the STERTO method  performs the best and our proposed algorithm has advantages over the FUSE method in terms of metrics R-SNR, CC, SAM, and ERGAS. The HySure method has comparable performance to our algorithm, but requires more computational time. Moreover, our proposed algorithm gets a  higher  R-SNR value and  lower SAM value when compared to the SCOTT. Figure \ref{fig3} also provides an intuitive and reasonable display of super-resolution images.

In addition,  we demonstrate the results given by  different algorithms for Indian pines in Table \ref{tab:2} under the conditions $d=6, q=9,$ and $d=4, q=5.$ It seems that the values $d$ and $q$ affect the performances of all methods in some degree. However, the ranking of each evaluation parameters of different methods  almost do not change when $d$ and $q$ vary.

\begin{table}
	\begin{tabular}{ll|l|ccccc}
		\toprule
		&&\multirow{2}*{Algorithm} & \multicolumn{5}{c}{quality evaluation metrics} \\ \cline{4-8}
		&& &\text { R-SNR }  & \text { CC } & \text { SAM } &\text { ERGAS }& \text { TIME(s) } \\
		\hline
		&&\text { best value } & $+\infty$  & 1 & 0 & 0 & - \\ \hline
        &&\text { STERTO  } & 17.1952 & 0.987 & 0.4548 & 4.3075 & 1.1343 \\
        &d=4&\text { SCOTT } & 18.8878 & 0.9903 & 0.3651 & 3.8203 & 0.1726 \\
		&&\text { HySure } &18.3815 &0.9890 &0.3519 &4.0037 & 9.0540\\
		&q=9&\text { FUSE} &9.5258 &0.8919 &0.3769 &12.2427& 0.1248 \\
		&&\text { TTDSR } & 17.1458 & 0.9858 & 0.1089 & 4.5147 & 1.5402 \\
		\hline
	\end{tabular}
\caption{Comparison of performance of different algorithms on Salinas-A scene.}\label{tab:3}
\end{table}
\begin{figure}[h]
	\centering
	\subfigure[\label{fig:a}]{
		\includegraphics[scale=0.68]{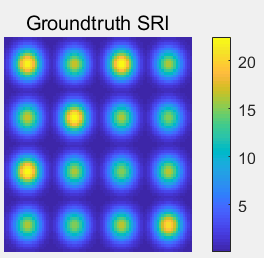}}
	\subfigure[\label{fig:c}]{
		\includegraphics[scale=0.68]{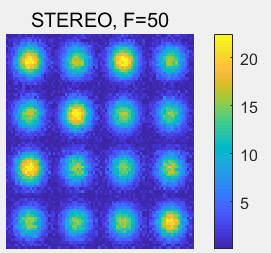}}
	\subfigure[\label{fig:e}]{
		\includegraphics[scale=0.68]{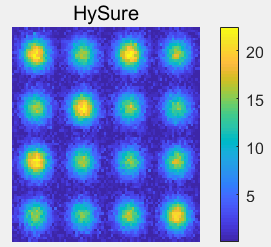}}
	\subfigure[\label{fig:b}]{
		\includegraphics[scale=0.68]{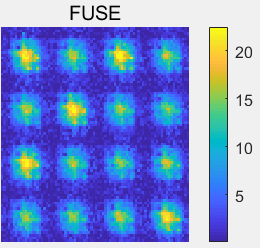}}
\subfigure[\label{fig:b}]{
		\includegraphics[scale=0.68]{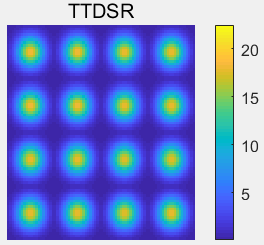}}
\caption{Spectral slice 120 of the SRI, Salinas-A scene.}
	\label{fig4}
\end{figure}

In the second example, Salinas-A scene dataset comes from the hyperspectral remote sensing data platform, which is available in \cite{28}. Similarly, white Gaussian noise  is added to $ \mathcal {Y}_ {h} $ and $ \mathcal {Y}_ {m}$ with an input SNR of $30$ dB. Consistently, we conduct experiments on the dataset under the conditions of $d=4 $, $q=9 $, and the results are shown in Table \ref{tab:3}. It is found that when $d=6, q=9 $ or $d=4, q=5 $ the metrics of the listed algorithms are almost consistent with the results of $d=4 $, $q=9. $  We omit results of under these two conditions. In Table \ref{tab:3}, compared with the STERTO, our proposed algorithm has comparable signal-to-noise ratio  and time. Compared to other algorithms, our method get the  lowest SAM value. Furthermore, it is evident that the TTDSR algorithm performs better than FUSE. In Figure \ref{fig4}, the super-resolution images obtained by different algorithms are shown. Due to the low SAM value, the image recovered by TTDSR algorithm are relatively clearer.
\subsection{Further numerical reuslts of TTDSR}
\begin{figure}[h]
	\centering
	\subfigure[\label{fig:a}]{
		\includegraphics[scale=0.40]{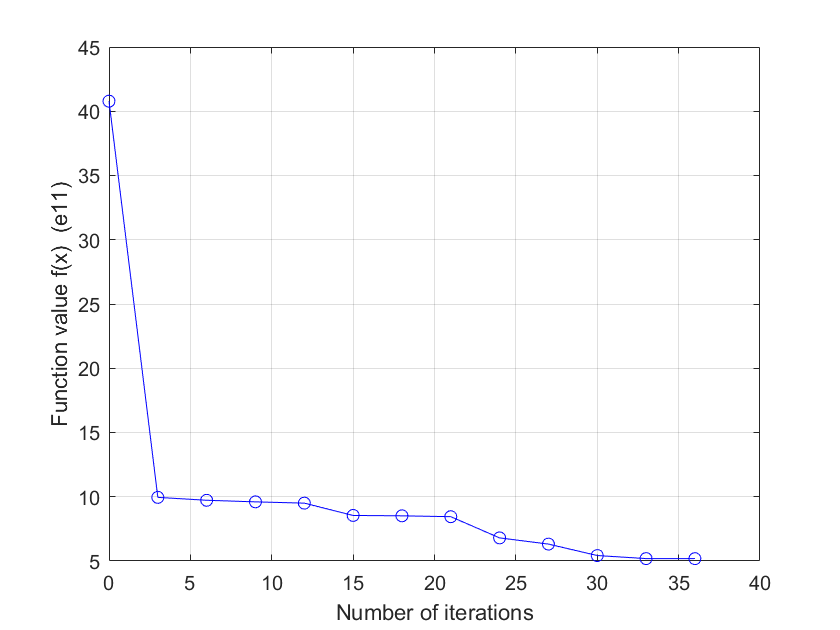}}
	\subfigure[\label{fig:b}]{
		\includegraphics[scale=0.40]{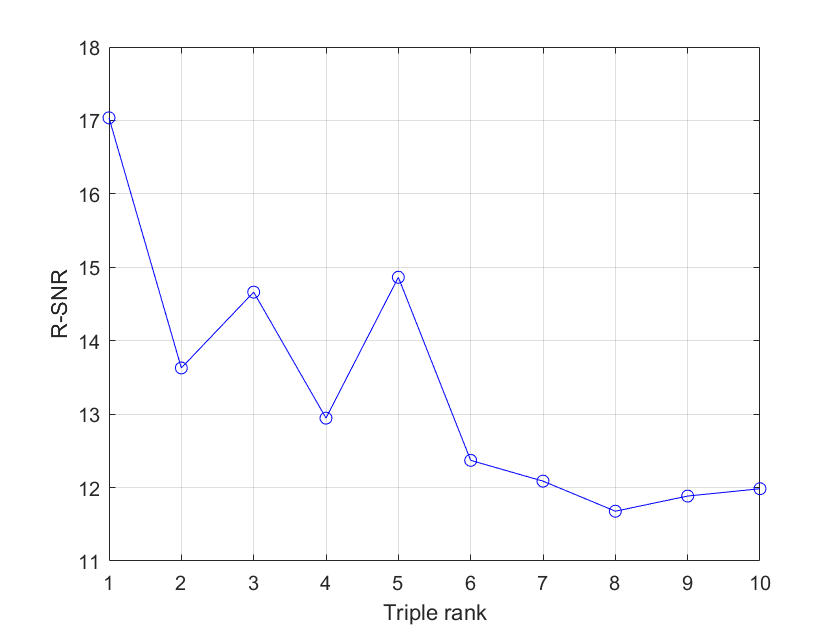}}
	\caption{Further experiments on Indian pines.}
	\label{fig5}
\end{figure}

In this subsection we further show the numerical results of TTDSR implemented on Indian pines dataset. The curve in Figure \ref{fig5}(a) displays the objective function value in each iteration, which verifies the theoretical conclusion that the sequence $\{f(\mathbf{x}_k)\}$ is decreasing.

In theory, the SRI is a low rank tensor. However, in the numerical experiments, we have no prior knowledge of  the triple rank of the SRI tensor and the rank is given artificially. For Indian pines dataset, we run TTDSR algorithm ten times with the triple rank changing from 1 to 10 accordingly. In Figure \ref{fig5}(b), we demonstrate the values of R-SNR corresponding to different triple ranks. In this example, the best R-SNR is attained when the triple rank is 1.

From the above analysis, we can see that among all algorithms, the HySure method performs the best but costs much more time than others. This is because it establishes an optimization problem of convex objective function with vector total variation regularization. The TV regularizer calculates the dispersion difference of the image in the horizontal and vertical directions. Our method has a significant advantage over the SCOTT in that we only need to consider a triple rank, while the rank of tucker decomposition \cite{28} is an array.
\section{Conclusion}
\label{section:7}
In this work, we provide a novel tensor triple decomposition model for hyperspectral super-resolution. Firstly, in order to capture the global interdependence between hyperspectral data of different modes, we use triple rank to characterize its low rank. Then we propose a optimization algorithm TTDSR to get the desired hyperspectral super-resolution image. Using the triple decomposition theorem, we cleverly obtain the gradient of the objective function of the model, which provides great help for solving the problem. Due to the algebraic nature of the objective function $f(\mathbf{x}) $, we apply the Kurdyka-$\L$ojasiewicz property in analyzing the convergence of the sequence of iterates generated by TTDSR.  In addition, experiments on two datasets show the feasibility and effectiveness of the TTDSR.  This work opens up a new prospect for realizing hyperspectral super-resolution by using various tensor decompositions.

\appendix
\section{The process of vectorization of the variable}
\label{AppendixA}
The gradient of the objective function in \eqref{3.3} is calculated by   vectorization of the variable as 	
\begin{equation}
\mathbf{x}:=\left(\begin{array}{l}
\mathbf{a}_{m_{1} r r} \\
\mathbf{b}_{m_{2} r r} \\
\mathbf{c}_{n_{3} r r}
\end{array}\right)=\left(\begin{array}{c}
A_{m_{1} r r \times 1} \\
B_{m_{2} r r \times 1} \\
C_{n_{3} r r\times 1}
\end{array}\right).
\end{equation}

We directly vectorize the optimal variables of (\ref{3.3}) in accordance with preset principles while dealing with small-scale data.
Firstly, noting that $ f\left(\mathbf{x}\right)= f_{1}\left(\mathbf{x}\right)+ f_{2}\left(\mathbf{x}\right)+ f_{3}\left(\mathbf{x}\right).$  Secondly, vectorizing the known tensors $\mathcal{Y}_{h}$, $ \mathcal{Y}_{m}$, we get
$$
\begin{aligned}
vec(\mathcal{Y}_{h})=\mathbf{d}_h \quad \quad  and \quad \quad vec( \mathcal{Y}_{m})=\mathbf{d}_m.
\end{aligned}
$$
The symbol $vec$ indicates the vectorization operator.
It yields from (\ref{a.1}) that
$$
\begin{aligned}
{\left(\mathbf{y}_h\right)}_{n_{1}n_{2}n_{3}}
&=\left(\left(\left(C_{lm \times n_{3}}\otimes E_{n_{2}}\right)^{\top}\otimes P_{n_{1}\times{m_{1}}}\right)\left(\left(E_{m} \otimes K_{n \times n_{2}l}\right)^{\top} \otimes E_{m_{1}}\right)\right) \mathbf{a}_{m_{1}nm} \\
&=\left(C_{n_{3} \times lm}\otimes E_{n_{2}}\otimes P_{n_{1}\times m_{1}}\right)\left(E_{m}\otimes K_{n_{2}l\times n}\otimes E_{m_{1}}\right)\mathbf{a}_{m_{1}nm}\\
&=\left(C_{n_{3} \times lm}\otimes P_{n_{1}n_{2}\times m_{1}n_{2}}\right)\left(E_{m} \otimes K_{n_{2}l\times n}\otimes E_{ m_{1}}\right)\mathbf{a}_{m_{1}nm}. \\
\end{aligned}
$$
Hence,
$$
f_{1}\left(\mathbf{x}\right)=\left\|\left(C_{n_{3} \times lm}\otimes P_{n_{1}n_{2}\times m_{1}n_{2}}\right)\left(E_{m} \otimes K_{n_{2}l\times n}\otimes E_{ m_{1}}\right)\mathbf{a}_{m_{1}nm}-\mathbf{d}_{h}\right\|_{2}^{2}.
$$
It is easy to see
$$
\begin{aligned}
\frac{\partial f_{1}(\mathbf{x})}{\partial \mathbf{a}_{m_{1}nm}}=2 &\left(\left(C_{n_{3} \times lm}\otimes P_{n_{1}n_{2}\times m_{1}n_{2}}\right)\left(E_{m} \otimes K_{n_{2}l\times n}\otimes E_{m_{1}}\right)\right)^{\top}\\
&\cdot\left(\left(C_{n_{3} \times lm}\otimes P_{n_{1}n_{2}\times m_{1}n_{2}}\right)\left(E_{m} \otimes K_{n_{2}l\times n}\otimes E_{ m_{1}}\right)\mathbf{a}_{m_{1}nm}-\mathbf{d}_{h}\right) \\
=& 2\left(E_{m} \otimes K_{n\times n_{2}l}\otimes E_{m_{1} }\right)\left(C_{ lm\times n_{3}}\otimes P_{m_{1}n_{2}\times n_{1}n_{2} }\right)\\
&\cdot\left(\left(C_{n_{3} \times lm}\otimes P_{n_{1}n_{2}\times m_{1}n_{2}}\right)\left(E_{m} \otimes K_{n_{2}l\times n}\otimes E_{ m_{1}}\right)\mathbf{a}_{m_{1}nm}-\mathbf{d}_{h}\right).
\end{aligned}
$$
Similarly,  it yields from (\ref{a.4}) that
$$
\begin{aligned}
{\left(\mathbf{y}_m\right)}_{m_{1}m_{2}m_{3}} &=\left(\left(\left(E_{m} \otimes B_{n \times  m_{2}l}\right)\left(G_{lm \times m_{3}} \otimes E_{m_{2}}\right)\right)^{\top} \otimes E_{m_{1}}\right) \mathbf{a}_{m_{1}nm}\\
&=\left(\left(G_{m_{3}\times lm} \otimes E_{m_{1}m_{2}}\right)\left(E_{m} \otimes B_{m_{2} l \times n}\otimes E_{m_{1}}\right)\right)\mathbf{a}_{m_{1}nm}.
\end{aligned}
$$
Hence,
$$
f_{2}\left(\mathbf{x}\right)=\left\|\left(\left(G_{m_{3}\times lm} \otimes E_{m_{1}m_{2}}\right)\left(E_{m} \otimes B_{m_{2}l \times n}\otimes E_{m_{1}}\right)\right)\mathbf{a}_{m_{1}nm}-\mathbf{d}_{m}\right\|_{2}^{2},
$$
and
$$
\begin{aligned}
\frac{\partial f_{2}(\mathbf{x})}{\partial \mathbf{a}_{m_{1}nm}}=2 &\left(\left(G_{m_{3}\times lm}\otimes E_{m_{1}m_{2}}\right)\left(E_{m} \otimes B_{ m_{2}l\times n}\otimes E_{m_{1}}\right)\right)^{\top}\\
&\cdot\left(\left(G_{m_{3}\times lm}\otimes E_{m_{1}m_{2}}\right)\left(E_{m} \otimes B_{m_{2}l\times n}\otimes E_{m_{1}}\right)\mathbf{a}_{m_{1}nm}-\mathbf{d}_{m}\right) \\
=2 &\left(E_{m}\otimes B_{n\times m_{2}l}\otimes E_{m_{1}}\right)\left(G_{lm\times m_{3}}\otimes E_{m_{1}m_{2}}\right)\\
&\cdot\left(\left(G_{m_{3}\times lm}\otimes E_{m_{1}m_{2}}\right)\left(E_{m} \otimes B_{m_{2}l\times n}\otimes E_{m_{1}}\right)\mathbf{a}_{m_{1}nm}-\mathbf{d}_{m}\right).
\end{aligned}
$$
For $
f_{3}\left(\mathbf{x}\right)=\mu\left\| \mathbf{x} \right\|_{2}^{2}$, we have
$$
\frac{\partial f_{3}(\mathbf{x})}{\partial \mathbf{a}_{m_{1} nm}}=2 \mu\mathbf{a}_{m_{1} nm}.
$$
In a word,
\begin{equation}
\begin{split}
\frac{\partial f(\mathbf{x})}{\partial \mathbf{a}_{m_{1} nm}}&= 2 \left(E_{m} \otimes K_{n\times n_{2}l}\otimes E_{m_{1}}\right)\left(C_{lm\times n_{3}}\otimes P_{m_{1}n_{2}\times n_{1}n_{2} }\right)\\
&\cdot\left(\left(\left(C_{n_{3}\times lm}\otimes P_{n_{1}n_{2}\times m_{1}n_{2}}\right)\left(E_{m} \otimes K_{n_{2}l\times n}\otimes E_{ m_{1}}\right)\right)\mathbf{a}_{m_{1}nm}-\mathbf{d}_{h}\right)\\
&+ 2 \left(E_{m}\otimes B_{n\times m_{2}l}\otimes E_{m_{1}}\right)\left(G_{lm\times m_{3}}\otimes E_{m_{1} m_{2}}\right)\\
&\cdot\left(\left(\left(G_{m_{3}\times lm}\otimes E_{m_{1}m_{2}}\right)\left(E_{m} \otimes B_{m_{2}l\times n}\otimes E_{m_{1}}\right)\right)\mathbf{a}_{m_{1}nm}-\mathbf{d}_{m}\right)\\
&+ 2 \mu\mathbf{a}_{m_{1} nm}.
\end{split}
\end{equation}
It yields from (\ref{a.2}) and (\ref{a.5}) that
$$
f_{1}\left(\mathbf{x}\right)=\left\|\left(\left(H_{n_{1}\times mn} \otimes P_{n_{2}n_{3}\times m_{2}n_{3}}\right)\left(E_{n} \otimes C_{n_{3}m \times l}\otimes E_{ m_{2}}\right)\right)\mathbf{b}_{m_{2}ln}-\mathbf{d}_{h}\right\|_{2}^{2},
$$
$$
f_{2}\left(\mathbf{x}\right)=\left\|\left(\left(A_{m_{1}\times mn} \otimes E_{m_{2}m_{3}}\right)\left(E_{n} \otimes G_{m_{3}m \times l}\otimes E_{m_{2}}\right)\right)\mathbf{b}_{m_{2}ln}-\mathbf{d}_{m}\right\|_{2}^{2},
$$
and
$$
f_{3}\left(\mathbf{x}\right)=\mu\left\|  \mathbf{x} \right\|_{2}^{2}.
$$
Therefore,
\begin{equation}
\begin{aligned}
\frac{\partial f(\mathbf{x})}{\partial \mathbf{b}_{m_{2}ln }}&= 2 \left(E_{n} \otimes C_{ l \times n_{3}m}\otimes E_{m_{2}}\right)\left(H_{mn\times n_{1} } \otimes P_{ m_{2}n_{3}\times n_{2}n_{3}}\right)\\
&\cdot\left(\left(\left(H_{n_{1}\times mn} \otimes P_{n_{2}n_{3}\times m_{2}n_{3}}\right)\left(E_{n} \otimes C_{ n_{3}m \times l}\otimes E_{ m_{2}}\right)\right)\mathbf{b}_{m_{2}LN}-\mathbf{d}_{h}\right)\\
&+ 2 \left(E_{n}\otimes G_{l\times m_{3}m }\otimes E_{m_{2}}\right)\left(A_{mn\times m_{1}}\otimes E_{ m_{2} m_{3}}\right)\\
&\cdot\left(\left(\left(A_{m_{1}\times mn}\otimes E_{m_{2}m_{3}}\right)\left(E_{n}\otimes G_{m_{3}m\times l }\otimes E_{m_{2}}\right)\right)\mathbf{b}_{m_{2}ln}-\mathbf{d}_{m}\right)\\
&+ 2 \mu\mathbf{b}_{m_{2} ln}.
\end{aligned}
\end{equation}
It also yields from (\ref{a.3}) and (\ref{a.6}) that
$$
f_{1}\left(\mathbf{x}\right)=\left\|\left(\left(H_{n_{1}\times nm} \otimes E_{n_{2}n_{3}}\right)\left(E_{m} \otimes K_{ n_{2}n \times l}\otimes E_{n_{3}}\right)\right)\mathbf{c}_{n_{3}lm}-\mathbf{d}_{h}\right\|_{2}^{2},
$$
$$
f_{2}\left(\mathbf{x}\right)=\left\|\left(\left(A_{m_{1}\times nm} \otimes P_{{m_{2}m_{3}}\times m_{2}n_{3}}\right)\left(E_{m} \otimes B_{m_{2} n\times l}\otimes E_{ n_{3}}\right)\right)\mathbf{c}_{n_{3}lm}-\mathbf{d}_{m}\right\|_{2}^{2},
$$
and
$$
f_{3}\left(\mathbf{x}\right)=\mu\left\| \mathbf{x} \right\|_{2}^{2}.
$$
Thus,
\begin{equation}
\begin{split}
\frac{\partial f(\mathbf{x})}{\partial \mathbf{c}_{n_{3} lm}}&= 2 \left(E_{m} \otimes K_{ l\times n_{2}n}\otimes E_{n_{3}}\right)\left(H_{nm\times n_{1} } \otimes E_{n_{2}\times n_{3}}\right)\\
&\cdot\left(\left(\left(H_{n_{1}\times nm} \otimes E_{n_{2}\times n_{3}}\right)\left(E_{m} \otimes K_{ n_{2}n \times l}\otimes E_{n_{3}}\right)\right)\mathbf{c}_{n_{3}lm}-\mathbf{d}_{h}\right)\\
&+ 2 \left(E_{m}\otimes B_{l\times m_{2}n}\otimes E_{n_{3}} \right)\left(A_{nm \times m_{1} }\otimes P_{m_{2}n_{3}\times m_{2}m_{3}}\right)\\
&\cdot\left(\left(\left(A_{m_{1}\times nm}\otimes P_{m_{2}m_{3}\times m_{2}n_{3}}\right)\left(E_{m}\otimes B_{ m_{2}n \times l}\otimes E_{n_{3}}\right)\right)\mathbf{c}_{n_{3}lm}-\mathbf{d}_{m}\right)\\
&+ 2 \mu\mathbf{c}_{n_{3}lm}.
\end{split}
\end{equation}
Then, its gradient is
\begin{equation}
\nabla f(\mathbf{x})= 2 \left(\begin{split}
&\left(E_{m} \otimes K_{n\times n_{2}l}\otimes E_{m_{1}}\right)\left(C_{lm\times n_{3}}\otimes P_{m_{1}n_{2}\times n_{1}n_{2} }\right)\\
&\cdot\left(\left(\left(C_{n_{3}\times lm}\otimes P_{n_{1}n_{2}\times m_{1}n_{2}}\right)\left(E_{m} \otimes K_{n_{2}l\times n}\otimes E_{ m_{1}}\right)\right)\mathbf{a}_{m_{1}nm}-\mathbf{d}_{h}\right)\\
&+ 2 \left(E_{m}\otimes B_{n\times m_{2}l}\otimes E_{m_{1}}\right)\left(G_{lm\times m_{3}}\otimes E_{ m_{1} m_{2}}\right)\\
&\cdot\left(\left(\left(G_{m_{3}\times lm}\otimes E_{m_{1}m_{2}}\right)\left(E_{m} \otimes B_{ m_{2}l\times n}\otimes E_{m_{1}}\right)\right)\mathbf{a}_{m_{1}nm}-\mathbf{d}_{m}\right)\\
&+ 2 \mu\mathbf{a}_{m_{1} nm}\\
&\rule[-10pt]{9cm}{0.05em}\\
&\left(E_{n} \otimes C_{ l \times n_{3}m}\otimes E_{m_{2}}\right)\left(H_{mn\times n_{1} } \otimes P_{ m_{2}n_{3}\times n_{2}n_{3}}\right)\\
&\cdot\left(\left(\left(H_{n_{1}\times mn} \otimes P_{n_{2}n_{3}\times m_{2}n_{3}}\right)\left(E_{n} \otimes C_{ n_{3}m \times l}\otimes E_{ m_{2}}\right)\right)\mathbf{b}_{m_{2}ln}-\mathbf{d}_{h}\right)\\
&+ 2 \left(E_{n}\otimes G_{l\times m_{3}m }\otimes E_{m_{2}}\right)\left(A_{mn\times m_{1}}\otimes E_{m_{2} m_{3}}\right)\\
&\cdot\left(\left(\left(A_{m_{1}\times mn}\otimes E_{m_{2}m_{3}}\right)\left(E_{n}\otimes G_{m_{3}m\times l  }\otimes E_{m_{2}}\right)\right)\mathbf{b}_{m_{2}ln}-\mathbf{d}_{m}\right)\\
&+ 2 \mu\mathbf{b}_{m_{2} ln}\\
&\rule[-10pt]{9cm}{0.05em}\\
&\left(E_{m} \otimes K_{ l\times n_{2}n}\otimes E_{n_{3}}\right)\left(H_{nm\times n_{1} } \otimes E_{n_{2}\times n_{3}}\right)\\
&\cdot\left(\left(\left(H_{n_{1}\times nm} \otimes E_{n_{2}\times n_{3}}\right)\left(E_{m} \otimes K_{ n_{2}n \times l}\otimes E_{n_{3}}\right)\right)\mathbf{c}_{n_{3}lm}-\mathbf{d}_{h}\right)\\
&+ 2 \left(E_{m}\otimes B_{l\times m_{2}n}\otimes E_{n_{3}} \right)\left(A_{nm \times m_{1} }\otimes P_{m_{2}n_{3}\times m_{2}m_{3}}\right)\\
&\cdot\left(\left(\left(A_{m_{1}\times nm}\otimes P_{m_{2}m_{3}\times m_{2}n_{3}}\right)\left(E_{m}\otimes B_{ m_{2}n \times l}\otimes E_{n_{3}}\right)\right)\mathbf{c}_{n_{3}lm}-\mathbf{d}_{m}\right)\\
&+ 2 \mu\mathbf{c}_{n_{3}lm}
\end{split}\right).
\end{equation}
\section{We introduce five lemmas for the proof of Lemma \ref{lem5.2}}
\label{AppendixB}
First, we consider the BFGS update $(\ref{4.1}),(\ref{4.2}),(\ref{4.3})$.
\begin{lemma}\label{lemB.1}
	Suppose that $ H_{k+1} $ is generated by $(\ref{4.3})$. Then, we have
	\begin{equation}
	\left\|H_{k+1}\right\| \leq\left\|H_{k}\right\|\left(1+\frac{2 MN}{\epsilon}\right)^{2}+\frac{N^2}{\epsilon}.\label{B.1}
	\end{equation}
\end{lemma}
\begin{proof}
	If $ \mathbf{y}_{k}^{\top} \mathbf{s}_{k}<{\epsilon} $, we get  $\rho_{k}=0 $ and $ H_{k+1}=H_{k}$. Hence, the inequality (\ref{B.1}) holds. Next, we consider the case $ \mathbf{y}_{k}^{\top} \mathbf{s}_{k} \geq \epsilon $. Obviously,  $\rho_{k} \leq \frac{1}{\epsilon}$. From Lemma \ref{lem5.1} and $\left\{\mathbf{x}_{k}\right\}$ is bouned, suppose
	there exists a positive number $ N $ such that $\left\|\mathbf{s}_{k}\right\| \leq N $ and we get
	\begin{equation}
	\left\|\mathbf{y}_{k}\right\| \leq 2 M .\label{B.2}
	\end{equation}
	Since
	\begin{equation}
	\left\|V_{k}\right\| \leq 1+\rho_{k}\left\|\mathbf{y}_{k}\right\|\left\|\mathbf{s}_{k}\right\| \leq 1+\frac{2MN}{\epsilon}\quad \text { and } \quad\left\|\rho_{k} \mathbf{s}_{k} \mathbf{s}_{k}^{\top}\right\| \leq \rho_{k}\left\|\mathbf{s}_{k}\right\|^{2} \leq \frac{N^2}{\epsilon},
	\end{equation}
	we have
	\begin{equation}
	\left\|H_{k+1}\right\| \leq\left\|H_{k}\right\|\left\|V_{k}\right\|^{2}+\left\|\rho_{k} \mathbf{s}_{k} \mathbf{s}_{k}^{\top}\right\| \leq\left\|H_{k}\right\|\left(1+\frac{2MN}{\epsilon}\right)^{2}+\frac{N^2}{\epsilon}.
	\end{equation}
	Hence, the inequality (\ref{B.1}) is valid.
\end{proof}
\begin{lemma} \label{lemB.2}
	Suppose that $ H_{k} $ is positive definite and $ H_{k+1}$ is generated by BFGS $(\ref{4.1}),(\ref{4.2}),(\ref{4.3})$.  Let $ \lambda_{\min }(H) $ be the smallest eigenvalue of a symmetric matrix  $H$. Then, we get $ H_{k+1}$ is positive definite and
	\begin{equation}
	\lambda_{\min }\left(H_{k+1}\right) \geq \frac{\epsilon}{\epsilon+4 M^{2}\left\|H_{k}\right\|} \lambda_{\min }\left(H_{k}\right).\label{B.5}
	\end{equation}
\end{lemma}
\begin{proof}
	For any unit vector $\mathbf{v}$, we have
	\begin{equation}
	\mathbf{v}^{\top} H_{k+1} \mathbf{v}=\left(\mathbf{v}-\rho_{k} \mathbf{s}_{k}^{\top} \mathbf{v y}_{k}\right)^{\top} H_{k}\left(\mathbf{v}-\rho_{k} \mathbf{s}_{k}^{\top} \mathbf{v y}_{k}\right)+\rho_{k}\left(\mathbf{s}_{k}^{\top} \mathbf{v}\right)^{2} .
	\end{equation}
	Let $ t \equiv \mathbf{s}_{k}^{\top} \mathbf{v} $ and
	\begin{equation}
	\phi(t) \equiv\left(\mathbf{v}-t \rho_{k} \mathbf{y}_{k}\right)^{\top} H_{k}\left(\mathbf{v}-t \rho_{k} \mathbf{y}_{k}\right)+\rho_{k} t^{2} .
	\end{equation}
	Because  $ H_{k} $ is positive definite, $ \phi(t) $ is convex and attaches its minimum at $ t_{*}=   \frac{\rho_{k} \mathbf{y}_{k}^{\top} H_{k} \mathbf{v}}{\rho_{k}+\rho_{k}^{2} \mathbf{y}_{k}^{\top} H_{k} \mathbf{y}_{k}}$ . Hence,
	\begin{equation}
	\begin{aligned}
	\mathbf{v}^{\top} H_{k+1} \mathbf{v}  &\geq \phi\left(t_{*}\right) \\
	&=\mathbf{v}^{\top} H_{k} \mathbf{v}-t_{*} \rho_{k} \mathbf{y}_{k}^{\top} H_{k} \mathbf{v} \\
	&=\frac{\rho_{k} \mathbf{v}^{\top} H_{k} \mathbf{v}+\rho_{k}^{2}\left(\mathbf{y}_{k}^{\top} H_{k} \mathbf{y}_{k} \mathbf{v}^{\top} H_{k} \mathbf{v}-\left(\mathbf{y}_{k}^{\top} H_{k} \mathbf{v}\right)^{2}\right)}{\rho_{k}+\rho_{k}^{2} \mathbf{y}_{k}^{\top} H_{k} \mathbf{y}_{k}}\\
	&\geq \frac{\mathbf{v}^{\top} H_{k} \mathbf{v}}{1+\rho_{k} \mathbf{y}_{k}^{\top} H_{k} \mathbf{y}_{k}}>0,
	\end{aligned}
	\end{equation}
	where the penultimate inequality holds because the Cauchy-Schwarz inequality is valid for the positive definite matrix norm  $\|\cdot\|_{H_{k}} $,
	\begin{equation}
	\left\|\mathbf{y}_{k}\right\|_{H_{k}}\|\mathbf{v}\|_{H_{k}} \geq \mathbf{y}_{k}^{\top} H_{k} \mathbf{v}.
	\end{equation}
	So $ H_{k+1}$ is also positive definite. From (\ref{B.2}), it is easy to verify that
	\begin{equation}
	1+\rho_{k} \mathbf{y}_{k}^{\top} H_{k} \mathbf{y}_{k} \leq 1+\frac{4 M^{2}\left\|H_{k}\right\|}{\epsilon}.
	\end{equation}
	Thus, we have
	\begin{equation}
	\mathbf{v}^{\top} H_{k+1} \mathbf{v} \geq \frac{\epsilon}{{\epsilon}+4 M^{2}\left\|H_{k}\right\|} \lambda _{\min }\left(H_{k}\right).
	\end{equation}
	Hence, we get the validation of (\ref{B.5}).
\end{proof}

Second, we turn to L-BFGS. Regardless of the selection of $ \gamma_{k} $ in (\ref{4.6}), we get the following lemma.

\begin{lemma}\label{lemB.3}
	Suppose that the parameter $\gamma_{k} $ takes Barzilai-Borwein steps $(\ref{4.6})$ . Then, we have
	\begin{equation}
	\frac{\epsilon}{4 M^{2}} \leq \gamma_{k} \leq \frac{N^2}{\epsilon}. \label{B.12}
	\end{equation}
\end{lemma}
\begin{proof}
	If $ \mathbf{y}_{k}^{\top} \mathbf{s}_{k}<\epsilon$, we get $ \gamma_{k}=1 $ which satisfies the bounds in (\ref{B.12}) obviously. Otherwise, we have
	\begin{equation}
	\epsilon \leq \mathbf{y}_{k}^{\top} \mathbf{s}_{k} \leq\left\|\mathbf{y}_{k}\right\|\left\|\mathbf{s}_{k}\right\|.
	\end{equation}
	Recalling (\ref{B.2}), we get
	\begin{equation}
	\frac{\epsilon}{N} \leq\left\|\mathbf{y}_{k}\right\| \leq 2 M \quad \text { and } \quad \frac{\epsilon}{2 M} \leq\left\|\mathbf{s}_{k}\right\| \leq N.
	\end{equation}
	Hence, we have
	\begin{equation}
	\frac{\epsilon}{4 M^{2}} \leq \frac{\mathbf{y}_{k}^{\top} \mathbf{s}_{k}}{\left\|\mathbf{y}_{k}\right\|^{2}} \leq \frac{\left\|\mathbf{y}_{k}\right\|\left\|\mathbf{s}_{k}\right\|}{\left\|\mathbf{y}_{k}\right\|^{2}}=\frac{\left\|\mathbf{s}_{k}\right\|}{\left\|\mathbf{y}_{k}\right\|}=\frac{\left\|\mathbf{s}_{k}\right\|^{2}}{\left\|\mathbf{y}_{k}\right\|\left\|\mathbf{s}_{k}\right\|} \leq \frac{\left\|\mathbf{s}_{k}\right\|^{2}}{\mathbf{y}_{k}^{\top} \mathbf{s}_{k}} \leq \frac{N^2}{\epsilon}.
	\end{equation}
	which means that  $\gamma_{k}^{\mathrm{BB} 1}$ and $\gamma_{k}^{\mathrm{BB} 2}$ satisfy (\ref{B.12}).
\end{proof}

Third, based on Lemmas \ref{lemB.1}, \ref{lemB.2}, \ref{lemB.3}, we obtain two lemmas as follows.
\begin{lemma} \label{lemB.4}
	Suppose that the approximation of a Hessian's inverse $ H_{k} $ is generated by L-BFGS $(\ref{4.4}), (\ref{4.5})$. Then, there exists a positive constant  $ C_{M} \geq 1 $ such that
	\begin{equation}
	\left\|H_{k}\right\| \leq C_{M}.\label{B.16}
	\end{equation}
\end{lemma}
\begin{proof}
	From Lemma \ref{lemB.3} and (\ref{4.5}), we have $ \left\|H_{k}^{(0)}\right\| \leq \frac{N^2}{\epsilon} $. Then, for (\ref{4.4})  and Lemma \ref{lemB.1}, we get
	\begin{equation}
	\begin{aligned}
	\left\|H_{k}\right\|& \leq\left\|H_{k-1}\right\|\left(1+\frac{2 MN}{\epsilon}\right)^{2}+\frac{N^2}{\epsilon} \\
	& \leq \cdots \\
	& \leq\left\|H_{k}^{(0)}\right\|\left(1+\frac{2MN}{\epsilon}\right)^{2 l}+\frac{N^2}{\epsilon}\sum_{m=1}^{l-1}\left(1+\frac{2 MN}{\epsilon}\right)^{2 m} \\
	& \leq \frac{N^2}{\epsilon}\sum_{m=1}^{l}\left(1+\frac{2MN}{\epsilon}\right)^{2 m} \equiv C_{M}.
	\end{aligned}
	\end{equation}
	Thus (\ref{B.16}) holds.
\end{proof}

\begin{lemma}\label{lemB.5}
	Suppose that the approximation of a Hessian's inverse  $H_{k}$  is generated by L-BFGS $(\ref{4.4}), (\ref{4.5})$. Then, there exists a constant  $ 0<C_{m}<1 $ such that
	\begin{equation}
	\lambda _{\min }\left(H_{k}\right) \geq C_{m}.
	\end{equation}
\end{lemma}	
\begin{proof}
	From Lemma \ref{lemB.3} and (\ref{4.5}), we have $ \lambda_{\min }\left(H_{k}^{(0)}\right) \geq \frac{\epsilon}{4 M^{2}} .$ Moreover, Lemma \ref{lemB.4} means   $\left\|H_{k-m}\right\| \leq C_{M} $ for all $ m=1, \ldots, l $. Hence, Lemma \ref{lemB.2} implies
	\begin{equation}
	\lambda_{\min }\left(H_{k-m+1}\right) \geq \frac{\epsilon}{{\epsilon}+4 M^{2} C_{M}} \lambda_{\min }\left(H_{k-m}\right).
	\end{equation}
	Then,  we obtain
	\begin{equation}
	\begin{aligned}
	\lambda_{\min }\left(H_{k}\right) &\geq \frac{\epsilon}{{\epsilon}+4 M^{2} C_{M}} \lambda_{\min }\left(H_{k-1}\right) \\
	& \geq \cdots \\
	& \geq\left(\frac{\epsilon}{{\epsilon}+4 M^{2} C_{M}}\right)^{l} \lambda_{\min }\left(H_{k}^{(0)}\right) \\
	& \geq \frac{\epsilon}{4 M^{2}}\left(\frac{\epsilon}{{\epsilon}+4 M^{2} C_{M}}\right)^{l} \equiv C_{m}.
	\end{aligned}
	\end{equation}
	Finally, the proof of Lemma \ref{lem5.2}  is straightforward from Lemmas \ref{lemB.4} and \ref{lemB.5}.
\end{proof}	

%\section*{Use of AI tools declaration}
%The authors declare they have not used Artificial Intelligence (AI) tools in the creation of this article.
%The article does not use Artificial Intelligence (AI) tools.
%The authors declare they have used Artificial Intelligence (AI) tools in the creation of this article.
%AI tools used:
%How were the AI tools used?
%Where in the article is the information located?

%The acknowledgments section should not be numbered.

%\section*{Acknowledgments}
%The authors would like to thank the editor and the associated
%reviewers for their constructive comments and suggestions that helped in improving
%the quality of this manuscript.

\medskip
%The information below will be filled in by AIMS production staff.
%Received May 2023; revised August 2023; early access xxxx 20xx.
\medskip

\end{document}